\documentclass[12pt,leqno]{article}
\usepackage{amsthm,amsfonts,amssymb,amsmath,eufrak,oldgerm}
\usepackage{epsfig}
\renewcommand\a{\alpha}
\newcommand\ha{\hat \alpha}
\renewcommand\b{\beta}
\newcommand\tb{\tilde\beta}

\newcommand\R{\mathbb R}

\newcommand\D{{${\mathcal D}_\eps$}}
\def\lx{\langle x \rangle}

\def\eps{\varepsilon }
\newcommand\kernel{\hbox{\rm Ker}}
\newcommand\trans{T}
\newcommand\errfn{\textrm{errfn}}
\newcommand\br{\begin{rem}}
\newcommand\er{\end{rem}}
\newcommand\bp{\begin{pmatrix}}
\newcommand\ep{\end{pmatrix}}
\newcommand\be{\begin{equation}}
\newcommand\ee{\end{equation}}
\newcommand\ba{\begin{equation}\begin{aligned}}
\newcommand\ea{\end{aligned}\end{equation}}

\newcommand{\CalB}{\mathcal{B}}
\newcommand{\CalE}{\mathcal{E}}
\newcommand{\CalF}{\mathcal{F}}
\newcommand{\CalC}{\mathcal{C}}

\newcommand{\CalK}{\mathcal{K}}
\newcommand{\CalL}{\mathcal{L}}
\newcommand{\CalN}{\mathcal{N}}
\newcommand{\CalO}{\mathcal{O}}

\newcommand{\CalT}{\mathcal{T}}

\newcommand{\RR}{{\mathbb R}}
\newcommand{\BbbR}{{\mathbb R}}

\newcommand{\ZZ}{{\mathbb Z}}

\newcommand{\Id}{{\rm Id }}
\newcommand{\e}{{\epsilon}}
\newcommand{\Range}{{\rm Range }}
\newcommand{\Span}{{\rm Span }}

\newcommand{\sgn}{\text{\rm sgn}}
\newtheorem{theo}{Theorem}[section]
\newtheorem{prop}[theo]{Proposition}
\newtheorem{cor}[theo]{Corollary}
\newtheorem{lem}[theo]{Lemma}
\newtheorem{defi}[theo]{Definition}

\newtheorem{rem}[theo]{Remark}
\newtheorem{rems}[theo]{Remarks}

\numberwithin{equation}{section}
\title{
Hopf bifurcation of viscous shock waves in
compressible gas- and magnetohydrodynamics 
}

\author{\sc \small Benjamin Texier\thanks{
Universit\'e Paris 7, Denis Diderot Institut de Mathématiques de Jussieu, UMR 7586 Case 7012 2, place Jussieu;
texier@math.jussieu.fr:
Research of B.T.  was partially supported
under NSF grant number DMS-0300487.
} and
Kevin Zumbrun\thanks{Indiana University, Bloomington, IN 47405;
kzumbrun@indiana.edu:
Research of K.Z. was partially supported
under NSF grant number DMS-0300487.
 }}
\begin{document}

\maketitle

\begin{abstract}
Extending our previous results for artificial viscosity systems,
we show, under suitable spectral hypotheses, 
that shock wave solutions of compressible Navier--Stokes (cNS)
and magnetohydrodynamics (MHD) equations undergo 
Hopf bifurcation to nearby time-periodic solutions.
The main new difficulty associated with physical viscosity and
the corresponding absence of parabolic smoothing is the need
to show that the difference between nonlinear and linearized
solution operators is quadratically small in $H^s$ for data in $H^s$.
We accomplish this by a novel energy estimate carried out
in Lagrangian coordinates;
interestingly, this estimate is false in Eulerian coordinates.
At the same time, we greatly sharpen and simplify the analysis
of the previous work.
%The rest of the argument combines pointwise Green function 
%estimates obtained by Mascia and Zumbrun and a streamlined 
%version of the bifurcation framework
%developed in the authors' previous work.
\end{abstract}

%\clearpage
\tableofcontents
%\clearpage
%%%%%%%%%%%%%%%%%
\bigbreak
\section{Introduction}

A well-known phenomenon in combustion is the appearance 
of ``galloping'', ``spinning'', and ``cellular'' instabilities
of traveling detonation fronts, apparently corresponding with
Hopf bifurcation of the background solution; 
see, e.g., \cite{KS, LyZ1, LyZ2, TZ1, TZ2} and references therein.
Mathematically, there are two distinct aspects of this phenomenon.
The first is to verify the spectral scenario associated with Hopf
bifurcation, consisting of a conjugate pair $\lambda_\pm(\eps)=
\gamma(\eps)+i\tau(\eps)$ of complex eigenvalues
of the linearized operator about the wave crossing the imaginary
axis from stable (negative real part) to unstable side as bifurcation
parameter $\epsilon$ varies from negative to positive values.
Here, $\eps$ measures variation in physical parameters such
as heat release, rate of reaction, or strength of the detonation.

The second issue is to show that this spectral scenario
indeed corresponds at the nonlinear level to Hopf bifurcation: i.e.,
apparition of nearby time-periodic solutions ${\bf u}^a$
of approximate period $2\pi/\tau(\eps)$, $\eps=\eps(a)$,
branching from the steady solution at $\eps(0)=0$.
What makes this nontrivial is the absence of a spectral gap
between zero and the essential spectrum of the linearized
operator about the wave, which circumstance prevents the application of 
standard bifurcation theorems for PDE, as in, e.g., \cite{C, He, MM, VT}.

The first issue has been studied for detonations by numerical \cite{KS},
formal asymptotic \cite{Er, FW, BMR},  and Evans function 
\cite{LyZ1, LyZ2} techniques.
The second issue has been studied in a much more general setting 
in \cite{TZ1, TZ2, SS}.
In particular, it was shown in \cite{TZ2} 
for smooth shock solutions of conservation
laws with artificial (strictly parabolic) viscosity that a spectral
Hopf scenario implies actual nonlinear Hopf bifurcation.
This result was recently sharpened by Sandstede and Scheel \cite{SS}
to include
also exponential localization of the perturbed solution
(the result of \cite{TZ2} asserts only $(1+|x|)^{-1}$ decay)
and exchange of spectral stability.
The results of \cite{SS} were obtained by quite different
``spatial dynamics'' techniques,
viewing the problem as a spatial elliptic boundary value problem on
the space of time-periodic functions, 
whereas the original results of \cite{TZ2} were obtained by 
``temporal dynamics'' techniques centered around
the temporal evolution operator.

Our purpose in the present paper is to extend the latter, nonlinear,
results to the physical (nonstrictly parabolic) systems of compressible 
gas- and magnetohydrodynamics (MHD) with either ideal
or ``real'', van der Waals-type, equation of state by a modification
of the original, temporal dynamics argument of \cite{TZ2}.
At the same time, we show that exponential localization may by an improvement
of the central cancellation estimate be obtained
by temporal as well as spatial techniques, recovering the sharpened estimates
of \cite{SS}; 
moreover, this same improvement greatly simplifies the entire analysis.
We expect, but have not verified that exchange of stability
may be obtained by our techniques as well.

In contrast to the detonation case, this investigation is not driven
by known physical phenomena, but by the common mathematical structure 
with detonation.
That is, in this case the mathematics suggests a possible physical 
phenomenon, and not the reverse.
It would be very interesting to look numerically or experimentally
for its appearance, particularly in the rich settings of 
MHD or phase-transitional gas dynamics 
with van der Waals-type equation of state,
where instability, similarly as in the detonation case (but
in contrast to gas dynamics), frequently occurs \cite{T}.

\br\label{Ioss}
\textup{
A third approach to bifurcation in the absence of a spectral gap,
besides those described in \cite{TZ2, SS},  may
be found in a recent work of Kunze and Schneider \cite{KuS}
in which they analyze pitchfork bifurcation in the absence of a
spectral gap using weighted-norm methods like those used to
study stability for dispersive systems \cite{PW}.
This approach has been used to treat stability of dispersive--diffusive
scalar undercompressive shocks \cite{Do}; however,
it does not appear to generalize to shock waves in the system case.
}
\er

\subsection{Equations and assumptions}\label{assumptions}
Consider a one-parameter 
family of standing viscous shock solutions 
\be\label{prof}
U(x,t)=\bar U^\eps(x),
\qquad \lim_{z\to \pm \infty} \bar U^\eps(z)=U_\pm^\eps
\quad \hbox{\rm (constant for fixed $\eps$)},
\ee 
of a smoothly-varying family of conservation laws 
\begin{equation} \label{sys}
U_t =\CalF(\epsilon, U):= (B(\eps, U)U_{x})_x- F(\epsilon, U)_x,
\qquad U\in \RR^n,
\end{equation}
with associated linearized operators
\be\label{Ldef}
L_\epsilon :=\partial \CalF/\partial U|_{U=\bar U^\epsilon}
= 
\partial_x B^\eps(x) \partial_x
-\partial_x A^\eps(x)  
\ee
$B^\eps(x)=B(\bar U^\eps(x))$,
$A^\eps(x):= F_U(\bar U^\eps(x), \eps)$.
%TODO: put back something later????  (OMIT, I think...)
Equations \eqref{sys} are typically shifts  $B(\eps, U)=B(U)$,
$F(\eps, U):= f(U)-s(\eps) U$
of a single equation 
$$
U_t = (B(U)U_x)_x - f(u)_x
$$
written in coordinates $\tilde x = x-s(\epsilon)t$ moving with 
traveling-wave solutions $U(x,t)=\bar U^\eps(x-s(\eps)t)$
of varying speeds $s(\epsilon)$. 
Profiles $\bar U^\eps$ satisfy the standing-wave ODE
\be\label{ode}
B(U)U'=F(\eps, U)- F(\eps, U^\eps_-).
\ee

Denote
\be\label{pmvalues}
B^\eps_\pm:=\lim_{z\to \pm \infty} B^\eps(z)=F_U(U^\eps_\pm, \eps),
\qquad
A^\eps_\pm:=\lim_{z\to \pm \infty} A^\eps(z)=F_U(U^\eps_\pm, \eps),
\ee
and $A(\eps, U)=F_u(\eps, U)$.
\bigskip

We make the following structural assumptions. 
%on the equations.
\medskip

(A1) \quad
$U=\left(\begin{array}{c} U_1 \\ U_2\end{array}\right)$,
\quad $B=\left(\begin{array}{cc} 0 & 0 \\ 0 & b \end{array}\right)$, \quad
$b$ nonsingular, 
where $U\in \BbbR^n$, $U_1\in  \BbbR^{n-r}$, 
$U_2\in \BbbR^r$, and $b\in \BbbR^{r\times r}$;
moreover, the $u$-coordinate $F_1(\eps, U)$ of $F$ is {\it linear in} $U$
(strong block structure).
\medskip\medskip

(A2)\quad There exists a smooth, block-diagonal, 
positive definite matrix $A^0(\eps, U)$
such that $A^0_{11} A_{11}$ is symmetric
and $A^0_{22}b$ is positive definite but not necessarily
symmetric (symmetric hyperbolic--parabolicity).
\medskip\medskip

To (A1)--(A2), we add the following more detailed hypotheses.
Here and elsewhere, $\sigma(M)$ denotes the spectrum of a matrix
or linear operator $M$.
%\cite{MaZ3, Z4, Z2}.
\medbreak

(H0) \quad  $F, B\in C^{k}$, $k\ge 5$.
%TODO: need 5 to directly apply \cite{MaZ3}.  What used in \cite{MaZ4}?

(H1) \quad  $\sigma(A^\eps_{11})$ 
%(necessarily real and semisimple) 
real, semisimple, 
nonzero, and constant multiplicity.

(H2) \quad $\sigma(A^\eps_\pm)$ real, simple, and nonzero.

(H3)\quad
$\R\sigma \left( i\xi_j dF(U_\pm) - 
\xi^2 B(U_\pm)\right) \leq -\frac{\theta |\xi|^2}{1+ |\xi|^2}$, 
$\theta>0$, for $\xi \in\BbbR$.

(H4)  \quad  Considered as connecting orbits of \eqref{ode}, 
$\bar U^\eps$ lie in an $\ell$-dimensional manifold, $\ell\ge 1$, 
of solutions \eqref{prof}, obtained
as a transversal intersection of the unstable manifold at $U^\eps_-$
and the stable manifold at $U_+^\eps$.
\medbreak
\medskip\medskip

\br\label{laxrmk}
\textup{
Conditions (H1)--(H2) imply that $U^\eps_\pm$ are nonhyperbolic rest points of
ODE \eqref{ode} expressed in terms of the $v$-coordinate, 
whence, by standard ODE theory,
\be\label{profdecay}
|\partial_x^\ell (\bar U^\eps-U^\eps_\pm)(x)|\le Ce^{-\eta|x|},
\qquad
0\le \ell\le k+1,
\ee
for $x\gtrless 0$, some $\eta$, $C>0$; 
in particular, $|(\bar U^\eps)'(x)| \le Ce^{-\eta|x|}$.
Condition (H4) implies in part that
$\bar U^\eps$ is either of standard {\it Lax type}, 
$\ell=1$, or nonclassical {\it overcompressive type},
$\ell>1$,
i.e., the hyperbolic convection matrices 
$A^\eps_\pm:=F_U(U_\pm, \eps)$ 
at $\pm\infty$ have, respectively, $p-\ell$ negative and
$n-p$ positive eigenvalues for $1\le \ell\le p \le n$.
In the Lax case, the ($\ell=1$)-dimensional manifold of solutions
\eqref{prof} consist of translates of $\Bar U^\eps$,
and $p$ is the characteristic family associated with the shock. 
For further discussion, see \cite{ZH, MaZ3, Z2}.
}
\er

\br\label{gensymm}
\textup{
Similarly as in \cite{Z2, Z3}, it is readily checked that the
conclusions of this paper hold also with (A1)--(A2) replaced
by the weakened assumptions
that there exist an invertible change of coordinates $U\to (U_1, W_2)$
for which $F_1$ is a linear function of $U_1$, $W_2$ with $\partial_{U_1}F_1$ 
symmetrizable and $B(U)\partial_x U= 
\bp 0 \\  \beta(U)\partial_x W_2\ep$ with
$\Re \sigma \big((\partial W_2/\partial U_2) \beta \big)\le -\theta<0$
and (H1) modified as in \cite{MaZ3} to account for nonvanishing of $B_{21}$,
replacing $A_{11}$ by $A_*:=A_{11}-B_{21}B_{22}^{-1}A_{12}$;
%enter from MaZ3...
see Remarks \ref{genen} and \ref{superfluous}.
These assumptions are used to obtain the $H^s$ 
energy estimates of Section \ref{energy}.
}
\er

Conditions (A1)--(A2) (or the alternatives described in Remark \ref{genen})
are a slightly strengthened version of the
corresponding hypotheses of \cite{MaZ4, Z2, Z3} for general systems 
with ``real'', or partially parabolic viscosity, the 
difference lying in the strengthened
block structure condition (A1): in particular, the assumed linearity
of the $U_1$ equation.  Conditions (H0)--(H4) are 
%slightly weaker than the corresponding assumptions in \cite{MaZ4, Z2, Z3}.
the same as those in \cite{MaZ4, Z2, Z3}.
The class of equations satisfying our assumptions, 
though not complete, is sufficiently
broad to include many models of physical interest,
in particular compressible Navier--Stokes equations and the 
equations of compressible magnetohydrodynamics (MHD), expressed
in Lagrangian coordinates, with either ideal
or ``real'' van der Waals-type equation of state; 
see Section \ref{applications} for further discussion.
The role of assumptions (A1)--(A2) in our analysis 
is discussed further in Appendix \ref{eulerian}.
%NO, outdated I guess...
%And, equations outside this class do not appear to have
%the properties needed for our bifurcation analysis; see Section \ref{eulerian}.
%We regard its identification as one of the main contributions of this paper.

A simple example is
the equations of isentropic gas dynamics written in Lagrangian coordinates:
\begin{equation}
\left\{\begin{array}{l}
 v_t -u_x=0,\\ u_t + p_x =((\nu/v) u_{x})_x,
\end{array}\right.
\label{NSeg}
\end{equation}
where $v>0$ denotes specific volume, $u$ fluid velocity, $p=p(v)$ pressure, 
and $\nu>0$ the coefficient of viscosity, with $p_v<0$
corresponding to hyperbolicity of the associated first-order system.
Denoting $U_1=v$, $U_2=u$, and considering a family of traveling-wave
solutions $U(x,t)=\bar U^\eps(x-s(\eps)t)$ with $p_v(v^\eps_\pm)<0$,
we obtain after the change of coordinates $x\to x-s(\eps)t$ a family of
stationary solutions $U(x,t)=\bar U^\eps(x)$ and equations
\begin{equation}
\left\{\begin{array}{l}
 v_t - s(\eps)v_x- u_x=0,\\ u_t-s(\eps)u_x + p_x =((\nu/v) u_{x})_x
\end{array}\right.
\label{iseg}
\end{equation}
satisfying (A1)--(A2), (H0)--(H4).

\br\label{euler}
\textup{
Equations \eqref{NSeg} written in Eulerian coordinates, or
\begin{equation}
\left\{\begin{array}{l}
 \rho_t + (\rho u)_x=0,\\ 
(\rho u)_t + (\rho u^2 + p)_x =\nu u_{xx},
\end{array}\right.
\label{Eulereg}
\end{equation}
where $\rho=1/v$ denotes density {\it does not}
satisfy (A1)--(A2), (H0)--(H4), violating the second
part of the strong block
structure assumption (A1).
(Here, $U_1=\rho$, $U_2=m=\rho u$ is momentum, so that the
first equation is still linear; however, the
parabolic term $u_{xx}=(m/\rho)_{xx}$ is not in diagonal form.)
%TODO: note for later discussion-
%to put things another way, it is $u$ and not $U_2=m$ that
%experiences parabolic smoothing.
Indeed, we make crucial use in the analysis of the
Lagrangian formulation, as discussed in Appendix \ref{eulerian}.
%Remark \ref{noEuler}.
}
\er

%TODO: put back? or delete...?
%In gas dynamics or MHD, assumptions (H1)--(H4) are 
%satisfied for Lax-type shocks moving with
%nonzero speed relative to background fluid velocity.

\subsection{Spectral criteria}\label{spectral}
%TODO: describe/introduce the Evans function.
%Recall stability condition of \cite{MaZ3, MaZ4, Z2}.
%Give Hopf condition.

As discussed in \cite{ZH, MaZ3, Z1, Z2, Z3}, the linearized
operators $L_\eps$ have no spectral gap, since essential spectrum
accumulates at $\lambda=0$ on the imaginary axis.
In this situation, standard stability and bifurcation criteria 
based on isolated spectra of $L_\eps$ are replaced by
generalized versions expressed in terms of an
associated {\it Evans function}.
The Evans function $D_\eps (\lambda)$, defined as a Wronskian of functions 
spanning the decaying manifolds of solutions of the eigenvalue equation
\be\label{evalue}
(L_\eps-\lambda)u=0
\ee
associated with $L_\eps$ at $x\to +\infty$ and $x\to -\infty$ is an
analytic function with domain containing $\{\Re \lambda \ge 0\}$,
whose zeroes away from the essential spectrum
correspond in location and multiplicity with eigenvalues
of $L_\eps$.
Its behavior is also closely linked with that of the resolvent kernel
of $L_\eps$, i.e., the Laplace transform with respect to time of the
Green function $G$; see, e.g., \cite{AGJ, GZ, ZH, ZS, Z1, Z2} 
for history and further details.

In \cite{MaZ3}
there was established the following {\it stability criterion}.

\begin{prop}[\cite{MaZ3}]\label{stabcrit}
Let $\bar U^\eps$, \eqref{sys}  
be a family of traveling-waves and systems 
satisfying assumptions (A1)--(A2), (H0)--(H4).
Then, $\bar U^\eps$ is linearly stable from
$L^1 \to L^p$ if and only if
there exist precisely $\ell$ zeroes ($\ell$ as in (H4)) 
of $D_\eps(\cdot)$ 
in the nonstable half-plane $\R \lambda \ge 0$, 
necessarily at the origin $\lambda=0$.
\end{prop}

We define an analogous {\it Hopf bifurcation criterion} as

\medbreak

\quad (\D) \quad
On a neighborhood of $\{\Re \lambda \ge 0\}\setminus \{0\}$,
the only zeroes of $D$ are (i) a zero of multiplicity $\ell$ at
$\lambda=0$, and (ii) a crossing conjugate pair of zeroes
$\lambda_\pm(\eps)=\gamma(\eps)+i\tau(\eps)$ with
$\gamma(0)=0$, $\partial_\eps \gamma(0)>0$, and $\tau(0)\ne 0$.
\medbreak

\begin{rems}\label{simple}
\textup{
1. In the simplest, Lax case, under (A1)--(A2), (H0)--(H4),
simplicity of the root $\lambda=0$,
condition \D(i), is equivalent to 
$\langle \pi^\eps, (\bar u^\eps)'\rangle=\pi^\eps\cdot (u^\eps_+- u^\eps_-)
\ne 0$ for
$\pi^\eps$ (constant) orthogonal to $S(A^\eps_-)\cup U(A^\eps_+)$ 
\cite{GZ, MaZ2, Z1, Z2}.
Under the normalization $\langle \pi^\eps, (\bar u^\eps)'\rangle=1$,
operator $\Pi_0^\eps f:= (\bar u^\eps)' \langle \pi^\eps, f\rangle$
plays the role of a ``generalized spectral projection'' onto 
$$
\kernel L(\eps)=\Span \{(\bar u^\eps)'\},
$$
and $\pi^\eps$ the role of a generalized left eigenfunction \cite{MaZ3, ZH}.
Note that $\pi^\eps$ lies outside the domain
of $\Pi_0^\eps$, a consequence of the absence of a spectral gap.
A similar, but more complicated condition holds in the overcompressive
case \cite{MaZ3, Z1, Z2}.
}

\textup{
2. By Proposition \ref{stabcrit}, condition \D 
with (A1)--(A2), (H0)--(H4) implies that $\bar U^\eps$
is linearly stable for $\eps <0$ and unstable for $\eps>0$;
that is, there is a transition from stability to instability
at $\eps=0$.
}
\end{rems}

\subsection{Results}\label{results}

We introduce the following notation, to be used throughout the paper.

\begin{defi}
\textup{
Let $B_2\subset B_1$ and $X_2\subset X_1$ 
denote the Banach spaces determined by norms
$\|U\|_{B_1}:=\|U\|_{H^1}$,
$\|\partial_x U\|_{B_2}:=\|\partial_x U\|_{B_1}+ \|U\|_{L^1}$
and
\ba\label{X1}
\|U\|_{X_1}&:=
\|e^{\eta \lx}U\|_{H^2},\\
%NOTE: only need 2, to get 1 below...
\|\partial_x U\|_{X_2}&:=\|\partial_x U\|_{X_1}
+
\|e^{2\eta \lx} U\|_{H^1},
%NOTE: controls \|e^{\eta \lx} U\|_{L^\infty}, all we need...
\ea
where $\eta>0$ and $\lx:=(1+|x|^2)^{1/2}$.
}
\end{defi}

Our main result is the following theorem establishing Hopf bifurcation
from the steady solution $\bar U^\eps$ at $\eps=0$ under
bifurcation criterion (\D).
%the assumptions of Sections \ref{assumptions} and \ref{spectral}.

\begin{theo}\label{main}
Let $\bar U^\eps$, \eqref{sys}  
be a family of traveling-waves and systems 
satisfying assumptions (A1)--(A2), (H0)--(H4), and (\D).
Then, for $r \ge 0$, $\eta>0$ sufficiently small and $C>0$ sufficiently large,
there is a $C^1$ function
$\epsilon(r)$, $\e(0)=0$, 
and a $C^1$ family of time-periodic solutions  ${\bf U}^{r}(x,t)$
of \eqref{sys} with $\epsilon=\epsilon(r)$, of period $T(r)$,
$T(\cdot)\in C^1$, $T(0)=2\pi/\tau(0)$, with 
\be\label{expdecaybd}
C^{-1}r \le \|{\bf U}^r-\bar U^\eps\|_{X_1}\le Cr
\ee
for all $t$.
For Lax shocks, 
up to fixed translations in $x$, $t$,
these are the only time-periodic solutions nearby in $X_1$
with period $T\in [T_0, T_1]$ for any fixed $0<T_0<T_1<+\infty$;
%Alternative:
%For overcompressive shocks, they are part of a $C^1$ $(\ell-1)$-parameter
%family of solutions that are likewise unique up to translation in $x$, $t$.
%If $U_+^\eps\ne U_-^\eps$, they are the only 
%nearby solutions of the more general form ${\bf U}^r(x-\sigma^r t,t)$
%with ${\bf U}^r(x,\cdot)$ periodic.
%OLDER
if $U_+^\eps\ne U_-^\eps$, they are the only 
nearby solutions of the more general form ${\bf U}^r(x-\sigma^r t,t)$
with ${\bf U}^r(x,\cdot)$ periodic.
For overcompressive shocks, they are part of a $C^1$ $(\ell-1)$-parameter
family of solutions that are likewise unique up to translation in $x$, $t$.
%
%OLDEST
%about the linearized eigenspace associated with $\lambda_\pm(\eps)$,
%unique up to translation in $t$ among nearby solutions originating in $\CalC$.
%For overcompressive shocks, they are part of a $C^1$ $\ell$-parameter
%family defined on a cone $\CalC$ (defined explicitly in \eqref{Cdef}) 
%about the linearized eigenspace associated with $\lambda_\pm(\eps)$,
%unique up to translation in $t$ among nearby solutions originating in $\CalC$.
\end{theo}

\begin{rem}\label{localization}
\textup{
Bound \eqref{expdecaybd}, by
Sobolev embedding, includes also the result of {exponential localization},
$|{\bf U}^r-\bar U^\eps |\le Cre^{-\eta|x|}$.
}
\end{rem}

\subsection{Analysis}\label{analysis}

Theorem \ref{main} is proved using the general bifurcation framework
established in \cite{TZ2} together with pointwise Green function bounds
established in \cite{MaZ3} for general hyperbolic--parabolic systems.  
This is for the most part straightforward, demonstrating the
power and flexibility of the frameworks set up in \cite{TZ2} and \cite{MaZ3}.
However, there is an interesting and apparently general difficulty
associated with the absence of parabolic smoothing, namely, the need to
show that nonlinear source terms, defined as the difference between
nonlinear solutions and solutions of the linearized equations, 
are quadratically small in the $X_1$ norm relative to the 
$X_1$ norm of initial data, the key point being to control the
high-derivative norm $\|\cdot\|_{H^s}$ despite apparent loss of derivatives.
We discuss this issue in detail in Section \ref{energy}.

An analogous issue arises in the nonlinear stability
theory in going from the strictly parabolic to the hyperbolic--parabolic
case, with the difficulty again to control an apparent loss of derivatives.
This is resolved in \cite{MaZ2, MaZ4} with an auxiliary nonlinear
energy estimate.
Here, as there, we obtain the needed derivative control
by auxiliary energy estimates (Prop. \ref{energyest}).
However, these are of a rather different type, being local
rather than global in time, and measuring variational rather than
time-asymptotic properties.
They are also somewhat more delicate, depending strongly on
the structure of the hyperbolic part of the equations.
In particular, they hold in Lagrangian but not Eulerian coordinates,
indicating the importance of nonlinear transport effects; see
the discussion of Appendix \ref{eulerian}.
%TODO: omit, obsolete:
%The resolution of these issues we regard
%as the main mathematical contribution of this paper.

A second new feature in the analysis is the incorporation
of exponentially weighted norms, yielding \eqref{expdecaybd}.
We accomplish this by an improved way of accounting cancellation,
which at the same time greatly simplifies the analysis.
The key estimate in \cite{TZ2} corresponds heuristically to showing that
the kernel
\be\label{canc}
\CalK_y(x,y):=
\sum_{j=0}^\infty K_y(x, jT; y)
\ee
of the $y$-derivative of formal inverse 
\be\label{formal}
(\Id-e^{LT})^{-1}=\sum_{j=0}^\infty e^{LjT},
\ee
$L:=\partial_x^2-a\partial_x$,
converges in $L^\infty(L^2(x); y)$ uniformly for $0<C^{-1}\le T\le C$, 
where 
$$
K(x,t;y):=c_0t^{-1/2}e^{(x-y-at)^2/4t}
$$
is a convected heat kernel, with $a> 0$,
with decay of 
$\CalK_y(x,y)=\CalK_y(x-y,0)
$
in $|x-y|$ determining the ultimate decay rate for $|{\bf U}^r-\bar U^\eps |$.

Since $\|K(x,jT;y)\|_{L^2(x)}\equiv Ct^{-1/4}$, convergence in
\eqref{canc} cannot be absolute, but must involve cancellation.
This was detected in \cite{TZ2} by using the defining
equation $K_t+aK_x-K_{xx}=0$ to write
\ba\label{explain}
\sum_{j=1}^\infty K_y(x, jT; y)&\sim
T^{-1}\int_T^{+\infty} K_y(x, t; y)dt\\
&=
-a^{-1}T^{-1}\int_T^{+\infty} K_t(x, t; y)dt
+a^{-1}T^{-1}\int_T^{+\infty} K_{yy}(x, t; y)dt\\
&=
a^{-1}T^{-1} K(x, T; y) +a^{-1}T^{-1}\int_T^{+\infty} K_{yy}(x, t; y)dt,\\
\ea
and observing that $\int_T^{+\infty} K_{yy}(x, t; y)dt$ is convergent and
$\le C(1+|x-y|)^{-1}$.

On the other hand, we could just as well have repeated the process to estimate
\ba\label{continued}
\int_T^{+\infty} K_{yy}(x, t; y)dt &=
a^{-1}T^{-1} K_y(x, T; y) +a^{-1}T^{-1}\int_T^{+\infty} K_{yyy}(x, t; y)dt,\\
\ea
with $|\int_T^{+\infty} K_{yyy}(x, t; y)dt|\le C(1+|x-y|)^{-3/2}.$
Continuing the process, we find that $\CalK_y$ decays to any polynomial order,
with asymptotic expansion
\be\label{asymptotic}
\CalK_y(x,y)= T^{-1}\sum_{j=0}^\infty a^{-j}\partial_y^j K(x,T;y),
\ee
each term of which exhibits exponential decay, but for which we have
shown only polynomial spatial bounds 
and not convergence to zero of 
the truncation error of successive finite sums.
Moreover, there is still the issue of the ``continuization error''
introduced at the first step of \eqref{explain}.

In the present analysis, we sidestep these issues, effectively
summing to all orders expansion \eqref{asymptotic}, by working
at the level of the inverse Laplace transform formula
\be\label{LIT}
K(x,y;t)=\oint_{\Gamma}e^{\lambda t} K_\lambda(x,y)d\lambda
\ee
used to obtain pointwise bounds for the actual, variable-coefficient
problem, where $K_\lambda(x,y)$ is the resolvent kernel, or Laplace
transform of $K(x,t;y)$, and $\Gamma$ is an appropriate sectorial contour
contained strictly in the resolvent set of $L:=\partial_x^2-a\partial_x$
and lying in the strictly negative half-plane $\Re \lambda\le -\eta_0<0$
for $|\lambda |\ge r$.

%hence bounded away from $\lambda=0$.
Contributions from $\lambda \ge r>0$ are negligible, and likewise
from $x<y$ (recall, $a>0$), and so we may
focus on the low-frequency regime $|\lambda |\le r$ and $x>y$, 
where $K_\lambda(x,y)= c(\lambda) e^{\mu(\lambda)(x-y)}$
with
\ba\label{exp}
c(\lambda)&=(c_0+\lambda c_1 + \cdots)\\
\mu(\lambda)&= \lambda \tilde \mu(\lambda),\qquad
\tilde \mu = -a^{-1} + d_1\lambda + \cdots,
\ea
and thus 
$
\partial_y K_\lambda=\lambda \tilde \mu(\lambda)K_\lambda.
$
The assumption that $\Gamma$ lie in the resolvent set implies
further that 
$
\Re \mu \le -\eta_0<0
$
for $\lambda \in \Gamma \cap \{|\lambda|\le r\}$.

Carrying out the sum \eqref{canc} exactly, therefore,
within Laplace inversion formula \eqref{LIT},
on low frequencies $|\lambda\le r|$ and $x>y$,
we obtain a contribution to $\CalK_y(x,y)$ of
\ba\label{exactsum}
\oint_{\Gamma \cap \{|\lambda|\le r\}}\sum_{j=0}^J e^{\lambda j T} \partial_y K_\lambda(x,y)d\lambda
&=
\oint_{\Gamma \cap \{|\lambda|\le r\}}\Big(\sum_{j=0}^J \lambda e^{\lambda j T}\Big)
 \tilde \mu(\lambda) K_\lambda(x,y)d\lambda\\
&=
\oint_{\Gamma \cap \{|\lambda|\le r\}} (1-e^{(J+1)T\lambda})
\Big(\frac{\lambda}{ 1- e^{\lambda T}}\Big)\\
&\qquad \times
 \tilde \mu(\lambda) K_\lambda(x,y)d\lambda\\
&\to
\oint_{\Gamma\cap \{|\lambda|\le r\}} \Big(\frac{\lambda}
{1- e^{\lambda T}}\Big)
 \tilde \mu(\lambda) K_\lambda(x,y)d\lambda\\
&\sim
a^{-1}T^{-1}\oint_{\Gamma\cap \{|\lambda|\le r\}} 
 |K_\lambda(x,y)|d\lambda\\
&\le 
Ce^{-\eta_0|x-y|}
\ea
for some $C$, $\eta>0$, by \eqref{exp} and $\Re \mu\le -\eta_0$.
Not only does this argument yield exponential spatial decay bounds,
but it turns out to be much easier than the previous argument
to generalize to the full, variable-coefficient problem,
thus streamlining and sharpening the analysis at the same time.

%TODO: 1. What about excited term?  Can we treat in the same way?
%(with $\errfn$ times $\lambda$?)  I think so, at least about
%the exact evaluation!  (However, we don't get $|x-y|$ decay,
%but only $e^{-\eta |x|}$ times constant term in $y$, i.e.,
%just time-asymptotic projector (!) But!  This is perfectly good
%as an operator from $L^1$ to exponentially decaying functions.
%And, from exp. decaying, it converges uniformly in $L^2$, perfect!
%
%Wow!  Getting simpler and simpler!!!
%YES! THIS WAY!

\br\label{relation}
\textup{
Accounting cancellation in this way on the whole solution
at once is much simpler than breaking the solution into
individual components and estimating each separately as
was done in \cite{TZ2}.
%2. Rmk. that this shows better the relation to spatial dynamics,
%since everything pulled back to boundary, like elliptic case somehow...
%(OR?)
Notice further that we have essentially performed the inversion \eqref{formal}
directly via the spectral resolution formula
\be\label{specinv}
(\Id-e^{LT})^{-1}
\partial_x =
\oint_\Gamma (1-e^{\lambda T})^{-1}\partial_y 
G_\lambda d\lambda,
\ee
which has an ``elliptic'' flavor perhaps somewhat analogous 
to the spatial dynamics point of view.
However, in the full problem,
there is a boundary contribution coming from the
$J\to \infty$ term in \eqref{exactsum} for certain
(stationary eigen-) modes, that is not reflected
in formula \eqref{specinv}; see Remark \ref{boundary}.
}
\er

\subsection{Discussion and open problems}\label{discuss}

In the strictly parabolic case \cite{TZ2}, 
the analysis of the one-dimensional case extended easily
to yield a corresponding multidimensional result 
for shocks propagating along a cylinder of finite cross-section
with artificial (Neumann or periodic) boundary conditions.
In the hyperbolic--parabolic case, however,
our reliance on Lagrangian coordinates in deriving $H^s$ energy
estimates limits us for the moment to one dimension; see Remark
\ref{barrier} and the discussion of Appendix \ref{eulerian}.
%Not clear how to go to multi-d.  Very interesting issue,
%since in all other respects we are able to do it (and
%quite easily in fact...).
%TODO: IDEA- get the energy estimates along particle paths,
%but do the rest in Eulerian???  (No, doesn't quite work, it
%seems... interesting... what to do???  Could there be a ``partially
%Lagrangian coordinate'' in $x_1$ direction only???  This seems to
%be the issue... OR????
%OTHER: What about FT of nonlinear equations?  could this reduce
%us to one-dim, where we could do a Lagrangian transf.???
%keep thinking...
It would be very interesting (both physically and mathematically)
to remove this restriction, perhaps by the incorporation of
nonautonomous effects in the linearized equations as in \cite{D}.

A further issue when and whether bifurcation
condition (\D) is actually satisfied in physically interesting situations.
%well be vacuous for ideal gas dynamics, 
%TODO: vaccuous is wrong word I think... -K
It may well be that (\D) does not occur at all for ideal gas dynamics, 
nor any other kind of instability.
%where shock waves are rather stable.
However, it seems quite possible that (\D) can occur 
in the richer setting of MHD,
where Lax-type shock waves are known to be sometimes unstable even
for an ideal gas equation of state \cite{T},
% (Trakhinin).
or for phase-transitional gas dynamics with 
a van der Waals-type equation of state.
Numerical investigation of (\D) across a range of shock
and detonation waves is an important direction for future investigation.
Numerical evaluation of the Evans function 
may carried out in a well-conditioned fashioned, as described, e.g.,
in \cite{Br1, Br2, BrZ, BDG, HuZ, BHR}.
%TODO: ALSO, combustion investiations with Greg and Jeff? (if done...)

The analysis of the present paper serves also as a stepping-stone
to the closely related but more complicated detonation case.
With suitable elaboration, the argument extends to that case
\cite{TZ3}, yielding Hopf bifurcation, or ``galloping instability''
of viscous strong detonation waves in one dimension.
As in the shock case, this result of course is subject
to verification of the spectral bifurcation hypothesis (\D).
For detonation solutions of reacting gas dynamics, galloping
is expected to occur, and Evans function results \cite{LyZ1, LyZ2}
show that instability if it occurs must be of a (possibly degenerate)
Hopf type spectral configuration.

Finally, we mention the problem of determining stability of 
the bifurcating time-periodic waves whose existence is
established in Theorem \ref{main}.
%TODO: add some discussion somewhere of this point???
%As discussed in Section \ref{stab1}, 
In the absence of a spectral gap,
%ENDTODO
our method of analysis does not directly yield stability
as in the finite-dimensional ODE case, but at
best partial information on the location of point spectrum
associated with oscillatory modes, with stability presumably
corresponding to the standard condition $d\eps/dr>0$;
see \cite{SS} for a corresponding result in the strictly parabolic case.
The hope is that one could combine such information
with an analysis like that carried out for stationary
waves in \cite{ZH, MaZ3},
adapted from the autonomous to the time-periodic setting: 
that is, a generalized Floquet analysis in the PDE setting
and in the absence of a spectral gap.
This seems a particularly exciting direction for further development
of the theory.

\medskip

{\bf Plan of the paper.}
We begin in Section \ref{linear} by recalling
the linearized bounds of \cite{MaZ3, Z2, TZ2},
and establishing the new bound described in \eqref{exactsum}.
In Section \ref{energy}, we carry out the 
nonlinear energy estimates controlling higher derivatives 
in nonlinear source terms.  
In Section \ref{bifurcation}, 
we prove the main theorem,
carrying out the bifurcation analysis following the framework of \cite{TZ2}.
Finally, in Section \ref{applications}, we show that our results 
apply to gas- and magnetohydrodynamics.

\medbreak
{\bf Acknowledgement.}
Thanks to Bj\"orn Sandstede
and Arnd Scheel for their interest in this work
and for stimulating discussions on spatial dynamics 
and bifurcation in the absence of a spectral gap.
B.T. thanks Indiana University
for their hospitality during the collaborative visit in
which the analysis was carried out.

\section{Linearized estimates}\label{linear}

Assuming (A1)--(A2), (H0)--(H4) (alternatively, the weakened
hypotheses of Remark \ref{genen}), let $L_\eps$ as in the introduction
be the linearized operators \eqref{Ldef}, and $\lambda_\pm (\eps)$ 
the crossing eigenvalues of (\D).

\begin{lem}[TZ2]\label{efns}
Associated with eigenvalues $\lambda_\pm(\eps)$ of $L_\eps$
are right and left eigenfunctions $\phi^\eps_\pm$ 
and $\tilde \phi^\eps_\pm \in C^k(x,\eps)$, $k\ge 2$,
exponentially decaying in up to $q$ derivatives as
$x\to \pm \infty$,  and $L_\eps$-invariant projection
\be\label{proj}
\Pi f:= \sum_{j=\pm} \phi^\eps_j(x)\langle \tilde \phi^\eps_j, f\rangle
\quad
\ee
onto the total (oscillatory) eigenspace $\Sigma^\eps:=\Span \{\phi^\eps_\pm\}$,
bounded from $L^q$ or $B_2$ to $W^{k,p}\cap X_2$ 
for any $1\le q,p\le \infty$.  Moreover,
\be\label{zeromassefn}
\phi^\eps_\pm= \partial_x \Phi^\eps_\pm,
\ee
with $\Phi^\eps\in C^{k+1}$ exponentially
decaying in up to $k+1$ derivatives as $x\to \pm \infty$.
\end{lem}

%TODO: check that this is enough for more complicated $X_2$ defined here...
\begin{proof}
From simplicity of $\lambda_\pm$, and the fact \cite{MaZ3} that
they are bounded away from the essential spectrum of $L_\eps$, 
we obtain either by standard 
spectral perturbation theory \cite{Kat} or by direct Evans-function
calculations \cite{GJ1, GJ2, MaZ3, Z2} that there exist $\lambda_\pm(\cdot) $,
$\phi^\eps_\pm(\cdot)\in L^2$ with the same smoothness $C^k(\eps)$, $k\ge 5$,
assumed on $F$.  The exponential decay properties in $x$ then follow
by standard asymptotic ODE theory; see, e.g., \cite{GZ, Z2, Z3}.
Finally, recall the observation of \cite{ZH} that, by divergence form
of $L_\eps$, we may integrate  $L_\eps \phi= \lambda \phi$ from
$x=-\infty$ to $x=+\infty$ to obtain 
$\lambda \int_{-\infty}^{+\infty}\phi(x)dx=0$, and thereby
(since $\lambda_\pm \ne 0$ by assumption)
\be\label{zmass}
\int_{-\infty}^{+\infty}\phi_\pm(x)dx=0,
\ee
from which we obtain by integration \eqref{zeromassefn} with the stated
properties of $\Phi_\pm$.
From \eqref{zeromassefn} and representation \eqref{proj}, we 
obtain by H\"older's inequality the stated bounds on projection $\Pi$.
\end{proof}

Defining $\tilde \Pi^\eps:=\Id -\Pi^\eps$, $\tilde \Sigma^\eps:= \Range \tilde \Pi^\eps$,
and $\tilde L_\eps:= L_\eps\tilde \Pi^\eps$, denote by
\be\label{kernel}
G(x,t;y):= e^{L_\eps t}\delta_y(x)
\ee
the Green kernel associated with the linearized
solution operator $e^{Lt}$ of the linearized evolution equations
$u_t=L_\eps u$, and
\be\label{transkernel}
\tilde G(x,t;y):= e^{\tilde L_\eps t}\tilde \Pi \delta_y(x)
\ee
the Green kernel associated with the transverse linearized
solution operator $e^{\tilde L_\eps t}\tilde \Pi$.
By direct computation,
$G = \CalO +  \tilde G$, where
\be \label{O}
\CalO(x,t;y):= e^{(\gamma(\eps)+i\tau(\eps)) t}\phi_+(x) \tilde \phi_+^t(y)
+
e^{(\gamma(\eps)-i\tau(\eps)) t}\phi_-(x) \tilde \phi_-^t(y).
\ee

\subsection{Short time estimates}\label{shorttime}

\begin{lem}\label{masspres}
For $0\le t\le T$, $1\le p\le \infty$, $\eta>0$,
%$r\ge 0$, 
and some $C=C(T)$,
\ba\label{BDuhamel}
\|e^{L_\eps t} f \|_{L^p},\,
\|e^{\tilde L_\eps t}\tilde \Pi  f \|_{L^p}
&\le C \|f\|_{L^p}.\\
\ea
\ba\label{XDuhamel}
\|e^{\eta\lx} \int e^{L_\eps t}\partial_x f \|_{L^p},\,
\|e^{\eta\lx}\int e^{\tilde L_\eps t}\tilde \Pi \partial_x f \|_{L^p}
&\le C \|e^{\eta\lx}f\|_{L^p},\\
\ea
\ba\label{varXDuhamel}
\|e^{\eta\lx} \int \partial_{\eps, t}e^{L_\eps t}\partial_x f \|_{L^p},\,
\|e^{\eta\lx}\int \partial_{\eps,t} e^{\tilde L_\eps t}\tilde \Pi \partial_x f \|_{L^p}
&\le C \|e^{\eta\lx}f\|_{W^{2^p}},\\
\ea
where $\lx:=(1+|x|^2)^{1/2}$.
\end{lem}

\begin{proof}
From standard $C^0$ semigroup bound
$|e^{L_\eps t}|_{L^p\to L^p} \le C$
and properties $e^{L_\eps t}=e^{L_\eps t}\Pi + e^{\tilde L_\eps t}\tilde \Pi$
and $\|\Pi f\|_{L^p}\le |f|_{L^p}$,
we obtain \eqref{BDuhamel}.
Likewise, we may obtain integrated bounds \eqref{XDuhamel}($r=0$)
using the divergence form of $L_\eps $,
by integrating the linearized equations with respect to $x$
to obtain linearized equations $U_t=\CalL_\eps  U$ for integrated variable
$$
U(x,t,\eps):=\int_{-\infty}^x u(z,t,\eps)dz,
\qquad
u(\cdot,t):=e^{L_\eps t}\partial_x f,
$$
with linearized operator $\CalL_\eps :=-A^\eps(x)\partial_x + \partial_x^2$
of the same parabolic form as $L_\eps $, then applying 
standard $C^0$ semigroup estimates
(alternatively, more detailed pointwise bounds as in \cite{MaZ3})
to bound 
$$
\|e^{\eta\lx}U(\cdot, t, \eps)\|_{L^p}=
\|e^{\eta\lx}e^{\CalL_\eps t}f\|_{L^p}\le C\|e^{\eta \lx} f\|_{L^p},
$$
the $e^{\tilde L t}\tilde \Pi$ bound then following by
relation $e^{L_\eps t}=e^{L_\eps t}\Pi + e^{\tilde L_\eps t}\tilde \Pi$
together with
$\|e^{\eta\lx}\int \Pi \partial_x f\|_{L^p}\le |f|_{L^p}$.
We obtain \eqref{XDuhamel} ($r>0$)
by the change of variables $V=Ue^{\eta \lx}$, which, since
$\alpha:= e^{\eta\lx}$ satisfies $|(d/dx)^k\alpha|\le C|\alpha|$
for $k\ge 0$, converts the linearized equations to an equation with the
same principal part plus lower-order terms with bounded coeffients.

Finally, \eqref{varXDuhamel} follows from $\partial_t e^{L(\eps)t}=
 e^{L(\eps)t}L(\eps)$ and the variational equation
$(\partial_t-L)\partial_\eps U= (\partial_\eps L)U$,
$U(t):=e^{L(\eps t}U_0$, together with
$\|L(\eps)U_0\|_{L^p}$, $\|L(\eps)U\|_{L^p}\le C\|U_0\|_{W^{2,p}}$.
\end{proof}

\subsection{Pointwise Green function bounds}\label{ptwise}

We now develop the key cancellation estimates analogous to \eqref{exactsum}
of the introduction, adapting the pointwise semigroup
methods of \cite{ZH, MaZ3, Z2} to the present case.
Our starting point is the inverse Laplace transform representation
\be \label{inverseLT2}
G(x,t;y) = {1\over 2\pi i} \text{\rm P.V.}\int_{\eta-i\infty}^{\eta + i\infty}
e^{\lambda t}G_\lambda (x,y)d \lambda,
\ee
$\eta>0$ sufficiently large, established in \cite{MaZ3}.

Deforming the contour using analyticity of $G_\lambda$ \cite{MaZ3}
across oscillatory eigenvalues $\lambda_\pm(\eps)$
we obtain $G=\tilde G+ O$, where $O$, defined in \eqref{O}, 
is the sum of the residues of the integrand at $\lambda_\pm$,
and, for $\nu, r>0$ sufficiently small,
\ba \label{tildeGrep}
\tilde G(x,t;y) &= 
{1\over 2\pi i} \oint_\gamma e^{\lambda t}G_\lambda (x,y)d \lambda\\
&\qquad +
{1\over 2\pi i} \text{\rm P.V.}
\Big(
\int_{-\nu-i\infty}^{-\nu -ri}
+
\int_{-\nu+ ri}^{-\nu + i\infty}
\Big)
e^{\lambda t}G_\lambda (x,y)d \lambda\\
&=: \tilde G^{I} + \tilde G^{II},
\ea
where $\gamma$ is the counterclockwise arc of 
circle $\partial B(0,r)$ connecting
$-\nu-ri$ and $-\nu+ri$, 
and $G^I$ as in \cite{Z3} is the low-frequency and $G^{II}$
the high-frequency component of $\tilde G$.
Define associated solution operators $\tilde S^I(t,\eps)$ and
$S^{II}(t,\eps)$ 
by
\be\label{soln}
S^\beta(t)f(x):=  \int_{-\infty}^{+\infty}
\tilde G^\beta(x,t;y)f(y)dy
\ee
and 
$$
S:=S^I+S^{II}.
$$ 

Supressing the parameter $\eps$,
denote by $a_j^\pm$, $r_j^\pm$, $l^\pm_j$ the eigenvalues and
associated right and left eigenvectors of 
$A^\eps_\pm= F_u(u^\eps_\pm, \eps)$.
Following \cite{MaZ3}, 
let $a^{*}_j(x)$, $j=1,\dots,J\le (n-r)$  denote the eigenvalues of 
$$
A_*:=A_{11}-B_{21}B_{22}^{-1}A_{12}=A_{11}.
$$
Let $L(X,Y)$ denote the space of bounded linear operators
from Banach space $X$ to $Y$, equipped with the usual operator norm
$|\cdot|_{L(X,Y)}$.

\begin{prop} \label{HF}
Under assumptions (A1)--(A2), (H0)--(H4), (\D) 
(alternatively, those of Remark \ref{genen}),
for $0\le q,r,s$, $0\le 2q+s\le 4$, and $\nu>0$,
\be\label{derivbds}
\partial_{\eps}^q
\partial_{t}^r
\partial_x^s \tilde G^{II}(x,t;y)=
\sum_{p\le (2q+ r)+ s}
\Big( 
\sum_{j=1}^J O(e^{-\nu t})\delta_{x-\bar a_j^* t}(-y) 
 + O(e^{-\nu (|x-y|+ t)})
\Big)
\partial_y^p.
\ee
\end{prop}

\begin{proof}
The case $(q=r=0)$ is immediate from the bounds of \cite{MaZ3}.
Derivatives with respect to $\eps$ may be converted using the
variational equation
\be\label{variational}
(\lambda- L)(\partial_\eps G_\lambda)= (\partial_\eps L) G_\lambda
\ee
into two spatial derivatives, thus extending to case $(r=0)$.
Likewise, time-derivatives, appearing as powers of $\lambda$
at the level of the resolvent kernel, have the effect of single spatial
derivatives on ``hyperbolic'' blocks, and double spatial derivatives
on ``parabolic'' blocks, in the notation of \cite{MaZ3}, yielding by
the same estimates as in case $(q=r=0)$ the asserted bounds.  This
completes the proof of the general case.
\end{proof}

\begin{cor}\label{HFcancel}
Under assumptions (A1)--(A2), (H0)--(H4), (\D) 
(alternatively, the weakened assumptions of Remark \ref{genen}),
$\sum_{j=0}^\infty S^{II}(jT)$ converges uniformly 
and absolutely in operator 
norm  $|\cdot|_{L(X_1,X_1)}$,
for $\eps$ sufficiently small and $T$ in any compact set,
for any $0<\eta<\nu$, to a limit that is $C^1$ in $(\eps, T)$
with respect to the operator norm $L(B_1, B_1)$.
\end{cor}

\begin{proof}
Straightforward computation using \eqref{derivbds}.
\end{proof}

\begin{prop}\label{LF}
Under assumptions (A1)--(A2), (H0)--(H4), (\D) 
(alternatively, the weakened assumptions of Remark \ref{genen}),
$ \sum_{j=0}^\infty S^{I}(jT)$ converges uniformly 
%TODO: put limit 0\le r\le k$ everywhere... (global fix).
(but conditionally)
in the operator norm $|\cdot|_{L(X_2,B_1)}$,
for $\eps$ sufficiently small and $T$ in any compact set bounded
away from the origin, with limiting kernel satisfying
for $0\le q, r, s$,  $0\le 2q+s\le 4$, $\nu>0$,
\be\label{kerbd}
\partial_\eps^q \partial_T^r \partial_x^r \sum_{j=0}^\infty \tilde G_y^{I}(x,jT;y)
= O(e^{-\nu|x-y|}+ e^{-\nu|x|}).
\ee
\end{prop}

\begin{proof}
(Convergence).
Since (by Sobolev embedding)
$\partial_x \|e^{2\eta \lx} \cdot \|_{L^\infty}\subset X_2$,
it is sufficient to show convergence in
$L( \partial_x \|e^{2\eta \lx} \cdot \|_{L^\infty}, L^2)$, or,
equivalently,  convergence in
$L( \|e^{2\eta \lx} \cdot \|_{L^\infty}, L^2)$
of the operator $ \sum_{j=0}^\infty S^I \partial_x$ with kernel
$$
\sum_{j=0}^J \tilde G_y^{I}(x,jT;y),
$$
where, by \eqref{tildeGrep},
\be \label{yrep}
\tilde G^I_y(x,jT;y) =
{1\over 2\pi i} \oint_\gamma e^{\lambda jT}\partial_y G_\lambda (x,y)d \lambda,
\ee
hence (summing under the integral as described in the introduction)
\ba \label{ysum}
\sum_{j=0}^J \tilde G_y^{I}(x,jT;y)&=
{1\over 2\pi i} \oint_\gamma 
\Big(\frac{1}{ 1- e^{\lambda T}}\Big)
\partial_y G_\lambda (x,y)d \lambda\\
&\qquad -
{1\over 2\pi i} \oint_\gamma 
e^{(J+1)T\lambda} 
\Big(\frac{1}{ 1- e^{\lambda T}}\Big)
\partial_y G_\lambda (x,y)d \lambda\\
&=: I + II.
\ea

Taking $y\le 0$ for definiteness,
recall from \cite{MaZ3, Z2} that, for $\lambda$ 
sufficiently small, $\tilde G_\lambda(x,y)$ 
can in the Lax or overcompressive case
be expanded analytically at $\lambda=0$ as a sum of ``excited'' terms 
\ba\label{excited}
&\lambda^{-1}\phi_m(x)l_k^{-\trans}e^{\lambda \tilde \mu_k^-(\lambda)y},
\qquad
\tilde \mu_k^-= -(a_j^-)^{-1} + O(\lambda),
\ea
where $\phi_j(x)=O(e^{-\nu|x|})$ are stationary modes of the linearized
equations; ``scattering terms''
\ba\label{scattering}
r_j^\pm l_k^{-\trans}e^{\lambda \tilde \mu_k^-(\lambda)y + 
\lambda \tilde \mu_j^\pm(\lambda)x},
\qquad
\tilde \mu_\beta^\pm= -(a_\beta^\pm)^{-1} + O(\lambda);
\ea
and error terms 
\ba\label{excitederr}
O(e^{-\nu(|x|+|y|)})
\ea
and
\ba\label{scatteringerr}
O(\lambda e^{\lambda \tilde \mu_k^-(\lambda)y + 
\lambda \tilde \mu_j^\pm(\lambda)x}),
\ea
with symmetric expansions for $y\ge 0$.
These bounds may be converted by the Riemann Saddlepoint 
($\sim$ Stationary Phase) estimates to the pointwise Green
function bounds of \cite{ZH, MaZ3, Z2}.

Taking $y$-derivatives, we find that 
$$
\Big(\frac{1}{ 1- e^{\lambda T}}\Big) \partial_y \tilde G_\lambda
=-T^{-1} (\lambda^{-1}+O(1))\partial_y \tilde G_\lambda
$$
expands as $(Ta_k^-)^{-1}$ times the same excited and scattering
terms \eqref{excited} and \eqref{scattering}, plus
error terms of the same order \eqref{excitederr}, \eqref{scatteringerr},
plus new pole terms of order
\be\label{newpole}
\lambda^{-1}O(e^{-\nu(|x|+|y|)}).
\ee
Thus, by the same Riemann saddlepoint estimates used to bound
$\tilde G^I(x,t;y)$ in \cite{ZH, MaZ3, Z2}, we find that the
contribution to term $II$ of all except the new terms \eqref{newpole}
satisfies exactly the same bounds as $\tilde G^I(x,(J+1)T;y)$,
expanding for $t=(J+1)T$ as the sum of excited terms
\be\label{excitedtwo}
(Ta_k^-)^{-1}
\phi_j(x) l_k^{-t}
\Big(\errfn\left(\frac{y+a_k^{-}t}{\sqrt{4\beta_k^{-}t}}\right)
-\errfn \left(\frac{y-a_k^{-}t}{\sqrt{4\beta_k^{-}t}}\right)\Big),
\ee
$\beta_k^\pm>0$, scattering terms
\begin{equation}\label{scatteringtwo}
\begin{aligned}
&\sum_{k=1}^n 
(1+t)^{-1/2}O(e^{-(x-y-a_k^{-} t)^2/Mt} e^{-\eta x^+}) \\
&\qquad +
\sum_{a_k^{-} > 0, \, a_j^{-} < 0} 
\chi_{\{ |a_k^{-} t|\ge |y| \}}
(1+t)^{-1/2} O(e^{-(x-a_j^{-}(t-|y/a_k^{-}|))^2/Mt}
e^{-\eta x^+}), \\
&\qquad+
\sum_{a_k^{-} > 0, \, a_j^{+}> 0} 
\chi_{\{ |a_k^{-} t|\ge |y| \}}
(1+t)^{-1/2} O(e^{-(x-a_j^{+} (t-|y/a_k^{-}|))^2/Mt}
e^{-\eta x^-})  \\
\end{aligned}
\end{equation}
bounded by convected heat kernels,
where $x^+$ (resp. $x^-$) denotes the positive (resp. negative) part of $x$,
and a neglible error
\be\label{negligible}
O(e^{-\nu(|x-y|+t)}).
\ee
Terms \eqref{newpole} by the same argument give a time-independent
contribution of
%TODO: mention the below? (comes to play in fact for detonation case...)
%might be worthwhile to remark, since seems a little strange this way.
%YES, I think I prefer to point out that such terms in fact don't
%occur, if accounted appropriately.
%NOTE: in fact, this cancels with pole from initial time zero
%contribution, if we were to carefully analyze... Alternatively,
%we could move contour FIRST to negative half-plane, then sum
%inside the integral, thus incurring no pole...
\be\label{newpoletwo}
\sum f_m(x)g_m(y)=
O(e^{-\nu(|x|+|y|)})
\ee
up to negligible error \eqref{negligible}.

Noting that all terms except \eqref{excitedtwo} and \eqref{newpoletwo}
decay in $L^2(x)$ at least at Gaussian rate $(1+t)^{-1/4}$, independent
of $y$, while kernels \eqref{excitedtwo} integrated against an
exponentially localized function converges exponentially to
$$
\phi_j(x) l_k^{-t}
= O(e^{-\nu |x|}),
$$
we find that $II$ converges as $J\to \infty$ to a kernel with
the claimed bound $O(e^{-\nu |x|}+e^{\nu|x-y|})$.
Term $I$ (independent of $J$) by the elementary bound
$|\partial_y G_\lambda(x,y)|\le Ce^{-\nu|x-y|}$ for $\lambda$
bounded away from essential spectrum of $L$, and the fact that
$\gamma$ is bounded away from both $\lambda=0$ and the essential
spectrum of $L$, satisfies trivially the bound $I = O(e^{-\eta|x-y|})$,
completing the proof of convergence while at the same time
establishing \eqref{kerbd} for $q, r, s=0$.

(Derivative bounds). The remaining bounds \eqref{kerbd}($q,r,s\not \equiv 0$)
follow easily, by the observation that the limiting summands in term
$II$ are independent of $T$ and satisfy the same bounds after 
$x$- or $\eps$- differentiation as before,
while $I$, by 
(either using the variational equations \eqref{variational} to convert
$\eps$-derivatives to two spatial derivatives as above,
or estimating directly as in \cite{TZ2})
$$
\partial_\eps^q\partial_x^s\partial_y \tilde G_\lambda(x,y)=O(e^{-\nu|x-y|})
$$
%NOTE:
%see, for example, the discussion of $\eps$-derivatives in \cite{TZ2}.
%Or, could convert into two spatial derivatives again using var'l equation
%as we did above, though slightly less efficient...
%TODO: check correct reference number Prop. \ref{Greenbounds}
%
and
$$
|\partial_T \Big(\frac{1}{ 1- e^{\lambda T}}\Big)|\le C
$$
for $\lambda\in \gamma$ 
likewise satisfies the same bounds after as before.
\end{proof}

\br\label{HZ}
\textup{
Note that, except for the excited term \eqref{excitedtwo},
the estimates on $II$ are rather crude compared to the detail
in \cite{MaZ3, Z2}, while the estimate for $I$ is trivial.
This represents a considerable simplification from \cite{TZ2}.
}
\er

\br\label{uc}
\textup{
In the undercompressive case, there appear in the expansion
of $\tilde G_\lambda$ additional terms of order
$\lambda^{-1}e^{-\nu(|x|+|y|)}$ and
$e^{\lambda \tilde \mu_j^\pm(\lambda)x}e^{-\nu|y|}$
that do not fit in our analysis.
As discussed in \cite{ZH, Z2}, this case is essentially different.
}
\er

\br\label{boundary}
\textup{
Note that there is indeed a contribution from the limiting
term $II$ at infinity, as mentioned in Remark \ref{relation}.
}
\er

\begin{cor}\label{cancelcor}
Under assumptions (A1)--(A2), (H0)--(H4), (\D) 
(alternatively, the weakened assumptions of Remark \ref{genen}),
$\sum_{j=0}^\infty S(jT)$ converges uniformly 
in operator norm $|\cdot|_{L(X_2,B_1)}$,
for $\eps$ sufficiently small and $T$ in any compact set bounded
away from the origin, to a limit that is bounded
in $|\cdot|_{L(X_2,X_1)\cap L(B_2,B_1)}$ and $C^1(\eps,T)$
in  $|\cdot|_{L(B_2,B_1)}$ on the subspace $X_2$.
\end{cor}
%TODO: change if necessary to agree with notation of \cite{TZ2}.

\begin{proof}
By Corollary \ref{HFcancel} and Proposition \ref{LF},
it is sufficient to establish the second statement for $S^{I}$
using bounds \eqref{kerbd}.  
Further,  since $\partial_x L^1\subset B_2$ and (by Sobolev embedding)
$\partial_x \|e^{2\eta \lx} \cdot \|_{L^\infty}\subset X_2$,
it is sufficient to show boundedness in 
$L(\partial_x L^1, L^2)$ and
$L( \partial_x \|e^{2\eta \lx} \cdot \|_{L^\infty}, X_1)$, or,
equivalently, boundedness in $L( L^1, L^2)$ 
of the operator $\partial_(\eps, T)^q \sum_{j=0}^\infty S^I \partial_x$,
with kernel
$$
\sum_{j=0}^\infty \tilde \partial_{\eps, T}^q G_y^{I}(x,jT;y),
$$
$q=0,1$, and boundedness in 
$L(  \|e^{2\eta \lx} \cdot \|_{L^\infty}, \|e^{\eta \lx} \cdot \|_{L^\infty})$
of the operators $ \partial_x^r \sum_{j=0}^\infty S^I \partial_x$,
$0\le r\le 3$ with kernels
$$
\partial_x^r\sum_{j=0}^\infty \tilde G_y^{I}(x,jT;y),
$$
both routine consequences of bound \eqref{kerbd}.
\end{proof}

From Corollary \ref{cancelcor} we obtain the following important conclusion.

\begin{prop}\label{inverse}
Under assumptions (A1)--(A2), (H0)--(H4), (\D) 
(alternatively, the weakened assumptions of Remark \ref{genen}),
for $\eps$ sufficiently small, and $0<T_0\le T\le T_1$,
$(\Id -S(\eps,T))$ has a right inverse 
$$
(\Id-S(\eps, T))^{-1}\in L(X_2,X_1)
$$
that, restricted to its domain $X_2$, is $C^1$ 
in $\eps$, $T$ in the $L(B_2,B_1)$ norm.
\end{prop}
%TODO: change if necessary to agree with notation of \cite{TZ2}.

\begin{proof}
The first assertion follows by a standard telescoping sum argument,
setting $ (\Id-S(\eps, T))^{-1}:= \sum_{j=0}^\infty S(\eps, T)^{j}$,
the second by the final assertion of Corollary \ref{cancelcor}.
See \cite{TZ2}, Section 2, for further details.
%TODO: check ref.
\end{proof}

\section{Nonlinear energy estimates}\label{energy}

We next carry out $H^s$-energy estimates on 
the perturbation equations
\be\label{trunc}
U_t-L_\eps U= Q_\eps (U,U_x)_x
\ee
of \eqref{sys} about $\bar U^\eps$, to be used
to control higher derivatives in the 
fixed point iteration used to carry out our Lyapunov--Schmidt reduction
and bifurcation analyses.

By standard energy estimates, we have $H^s$ well-posedness,
$$
\|U(\cdot, T)\|_{H^s}\le C\|U_0\|_{H^s}
$$
for $T$ uniformly bounded.
Likewise, we have a formally quadratic linearized truncation error
$|Q|=O(|U|(|U|+|U_x|)$ for $|U|\le C$.
Our goal in this section, and what is far from evident in the
absence of parabolic smoothing, is
to establish a quadratic bound on the {\it linearization error}:
\be\label{quaderror}
\|U(\cdot, T)-e^{L_\eps T}U_0\|_{H^s}\le C\|U_0\|_{H^s}^2.
\ee
The corresponding bound does not hold for quasilinear hyperbolic
equations, nor as discussed in Appendix \ref{eulerian} for systems
of general hyperbolic--parabolic type, due to loss of derivatives.
However, it follows easily for systems satisfying
assumptions (A1)--(A2).

\subsection{${\rm H}^s$ linearization error}\label{key}

By the strong block structure assumption (A1), we may write \eqref{trunc}
more precisely as
\be\label{blocktrunc}
U_t + (A^\eps U)_x -
(B^\eps U_x)_x
= \partial_x \bp 0 \\ q^\eps_1(U)(U,U) +  q_2^\eps(U)(U, \partial_x U_2) \ep,
\ee
where
$
B^\eps=\bp 0 & 0\\
0 & b^\eps \ep 
$
and $q_j(\cdot)$ are $C^{k-2}$ bilinear forms.
Appealing to (A2), we may differentiate $\ell$ times and multiply by 
$$
A^{0,\eps}:=A^0(\Bar U^\eps)=\bp A^{0,\eps}_{11} & 0\\
0 & A^{0,\eps}_{22}\ep
$$
to obtain
\be\label{symmtrunc}
A^{0,\eps} \partial_x^\ell U_t + \tilde A^\eps \partial_x^\ell U_x + 
-(\tilde B^\eps \partial_x^{\ell+1} U)_x
= \Theta_\ell +  A^{0,\eps}\partial_x\bp 0 \\ \tilde q^\eps_\ell \ep,
\ee
where $A^{0,\eps}$ is symmetric positive definite, $\tilde A^\eps_{11}$
is symmetric, $\tilde B^\eps=\bp 0 & 0\\0 & \tilde b^\eps\ep$
with $\tilde b^\eps$ symmetric positive definite, and
\ba\label{tildeq}
\|\Theta_\ell \|_{L^2}&= O(\|U\|_{H^\ell}),\\
\tilde q_\ell&=
\partial_x^{\ell}
\big(q^\eps_1(U)(U,U) + q_2^\eps(U)(U, \partial_x U_2) \big).
\ea

\begin{lem}\label{moserbd}
For $\ell \ge 1$,
$\|\tilde q_\ell\|_{H^\ell}\le  
C\big(\|U\|_{H^\ell} + \|U_2\|_{H^{\ell+1}}\big)
\big( \|U\|_{H^\ell}+ \|U\|_{H^\ell}^\ell\big).
$
\end{lem}

\begin{proof} Standard application of 
Moser's inequality; see, e.g., \cite{Ta, Z3}.
\end{proof}

\begin{prop}\label{energyest}
Assuming (A1), (A2), (H0),  for 
$1 \le s \le k-1$, $0\le T\le T_0$ uniformly bounded, 
some $C=C(T_0)>0$,
and $U$ satisfying \eqref{trunc} with
%TODO: check needs in proof...
initial data $U(\cdot, 0)=U_0$
sufficiently small in $H^s$, 
\ba\label{en}
\|U(\cdot, T)\|_{H^s}^2+\int_0^T \|U_2(\cdot, t)\|_{H^{s+1}}^2\, dt
&\le C\|U_0\|_{H^s}^2,\\
\ea
%and
\ba\label{varen}
\|U(\cdot, T)-e^{LT}U_0\|_{H^s}&\le C\|U_0\|_{H^s}^2.
\ea
\end{prop}

\begin{proof}
By symmetric positive definiteness of $A^{0,\eps}$,
$$
\CalE(U):=(1/2)\sum_{\ell=1}^s 
\langle \partial_x^\ell U, A^{0,\eps}\partial_x^\ell U\rangle
$$
defines a norm equivalent to $\|\cdot\|_{H^s}$, 
i.e., $\CalE(\cdot)^{1/2}\sim \|\cdot\|_{H^s}$.
Applying \eqref{symmtrunc}, we find that
\ba\label{encalc}
\partial_t \CalE(U)&=
-\langle \partial_x^s U, \tilde A^\eps \partial_x^{s+1}U\rangle
+\langle \partial_x^s U, \partial_x(\tilde B^\eps \partial_x^s U)
\rangle \\
&\quad + 
\langle \partial_x^s U_2, A^{0,\eps}_{22}\partial_x \tilde q_s \rangle
+O(\|U\|_{H^s}^2)\\
&= 
\langle \partial_x^s U, (1/2)\tilde A^\eps_x \partial_x^{s}U\rangle
+\langle \partial_x^{s+1} U_2, 
\tilde b^\eps \partial_x^{s+1} U_2 \rangle \\
&\quad
- \langle \partial_x^{s+1} U_2, A^{0,\eps}_{22}\tilde q_s \rangle
- \langle \partial_x^{s} U_2, \partial_x A^{0,\eps}_{22}\tilde q_s \rangle
+O(\|U\|_{H^s}^2)\\
&\le -\theta \|U_2\|_{H^{s+1}}^2 +O(\|U\|_{H^s}^2)\\
&\quad +O(\|U_2\|_{H^{s+1}}) \big(\|U\|_{H^s} + \|U_2\|_{H^{s+1}}\big)
\big( \|U\|_{H^s}+ \|U\|_{H^s}^s\big)
\ea
for some $\theta>0$.  So long as $\|U\|_{H^s}$ remains sufficiently
small, this gives
\ba
\partial_t \CalE(U) &\le -(\theta/2) \|U_2\|_{H^{s+1}}^2
+O(\|U\|_{H^s}^2)\\
&\le -(\theta/2) \|U_2\|_{H^{s+1}}^2
+C\CalE,
\ea
from which \eqref{en} follows by Gronwall's inequality, in the form
$$
\CalE(U(T))+ (\theta/2)\int_0^T \|U_2\|_{H^{s+1}}^2(t) dt \le C_2\CalE(U_0).
$$

To obtain \eqref{varen}, observe that
$V(\cdot, t):=U(\cdot, t)- e^{L_\eps t}U_0$ satisfies
\be\label{vareq}
E_t + (A^\eps E)_x -
(B^\eps E_x)_x
= \partial_x \bp 0 \\ q^\eps_1(U)(U,U) +  q_2^\eps(U)(U, \partial_x U_2) \ep,
\ee
with $E(\cdot, 0)=0$, and thus
\be\label{Esymmtrunc}
A^{0,\eps} \partial_x^\ell E_t + \tilde A^\eps \partial_x^\ell E_x + 
-(\tilde B^\eps \partial_x^{\ell+1} E)_x
= \Theta_\ell^E +  A^{0,\eps}\partial_x\bp 0 \\ \tilde q^\eps_\ell \ep,
\ee
similarly as in \eqref{symmtrunc}, with 
$\|\Theta_\ell^E\|_{L^2}=O(\|E\|_{H^\ell})$.
Thus, calculating as in \eqref{encalc}, we obtain
\ba\label{Eencalc}
\partial_t \CalE(E)&=
-\langle \partial_x^s E, \tilde A^\eps \partial_x^{s+1}E\rangle
+\langle \partial_x^s E, \partial_x(\tilde B^\eps \partial_x^s E)
\rangle \\
&\quad + 
\langle \partial_x^s E_2, A^{0,\eps}_{22}\partial_x \tilde q_s \rangle
+O(\|E\|_{H^s} (\|E\|_{H^s}+ \|U\|_{H^s}^2)\\
&= 
\langle \partial_x^s E, (1/2)\tilde A^\eps_x \partial_x^{s}E\rangle
+\langle \partial_x^{s+1} E_2, 
\tilde b^\eps \partial_x^{s+1} E_2 \rangle \\
&\quad
- \langle \partial_x^{s+1} E_2, A^{0,\eps}_{22}\tilde q_s \rangle
- \langle \partial_x^{s} U_2, \partial_x A^{0,\eps}_{22}\tilde q_s \rangle
+O(\|E\|_{H^s} (\|E\|_{H^s}+ \|U\|_{H^s}^2)\\
&\le -\theta \|E_2\|_{H^{s+1}}^2 +
+O(\|E\|_{H^s} (\|E\|_{H^s}+ \|U\|_{H^s}^2)\\
&\quad +O(\|E_2\|_{H^{s+1}}) \big(\|U\|_{H^s} + \|U_2\|_{H^{s+1}}\big)
\big( \|U\|_{H^s}+ \|U\|_{H^s}^s\big),
\ea
which, so long as $\|U\|_{H^s}$ remains sufficiently small, gives
$$
\partial_t \CalE(E) \le -\theta \|E_2\|_{H^{s+1}}^2 
+C\|E\|_{H^s}^2 
+C\|U\|_{H^s}^2\Big( \|U\|_{H^s}^2 + \|U_2\|_{H^{s+1}}^2\Big)
$$
and thus, using $\|U\|_{H^s}(t)\le C\|U_0\|_{H^s}$ by \eqref{en},
$$
\partial_t \CalE(E) \le 
+C\CalE(E)
+ C_2\|U_0\|_{H^s}^2\big( \|U\|_{H^s}^2 + \|U_2\|_{H^{s+1}}^2\big).
$$
Applying Gronwall's inequality, and using $E(0)=0$, we thus obtain
$$
\CalE(E(T))\le C_3\|U_0\|_{H^s}^2
\int_0^T \big( \|U\|_{H^s}^2 + \|U_2\|_{H^{s+1}}^2\big)(t)dt,
$$
yielding the result by \eqref{en} and $\CalE(E)\sim \|E\|_{H^s}^2$.
\end{proof}

\br\label{easy}
\textup{
Note that the above, finite-time estimate is considerably
simpler than the ``Kawashima-type'' global-in-time estimates 
used in the nonlinear stability analysis \cite{MaZ4, Z2, Z3},
which require also additional assumptions 
of, among other things, symmetrizability and ``genuine noncoupling'' of 
$A^\eps_\pm $ and $B^\eps_\pm$.
}
\er

\subsection{Weighted variational bounds}\label{lipbounds}

Define now
\be\label{N}
\CalN(U_0, \eps, T):=
U(\cdot, T) - e^{\tilde L(\eps)T}U_0.
\ee
where $U$ satisfies \eqref{trunc} with initial data $U_0$.

\begin{prop}\label{Nest}
Assuming (A1), (A2), (H0),  for 
$1 \le s \le k-1$, $0\le T\le T_0$ uniformly bounded, 
some $C=C(T_0)>0$, $\CalN$ is uniformly bounded in $X_2$
and $C^1$ in $B_2$
with respect to $U_0$, $T$, $\eps$ 
for $\|U_0\|_{X_1}$, $|\eps|$ sufficiently small, with
\ba\label{Nvaren}
\|\CalN(U_0, \eps, T)\|_{X_2}
&\le C\|U_0\|_{X_1}^2,\\
\|\partial_{U_0}\CalN(U_0, \eps, T)\|_{X_2}
&\le C\|U_0\|_{X_1},\\
\|\partial_{T,\eps}\CalN(U_0, \eps, T)\|_{B_2}
&\le C\|U_0\|_{X_1}^2.\\
\ea
\end{prop}

\begin{proof}
Similarly as in the proof of \eqref{XDuhamel}($r>0$),
we may make the change of variables 
$V=U\alpha(x)$, $\alpha(x):=e^{\eta \lx}$, $\lx:=(1+|x|^2)^{1/2}$,
to convert both \eqref{trunc} and the corresponding
linearized equations to equations for which the
same principal part plus lower-order terms with bounded coeffients
and the nonlinear part is $\alpha^{-1}$ times the same principal
part plus bounded factors times lower-order terms, thus recovering
$\|U(\cdot, T)-e^{LT}U_0\|_{X_1}\le C\|U_0\|_{X_1}^2$, or
\be\label{weightvar}
\|\CalN(U_0, \eps, T)\|_{X_1} \le C\|U_0\|_{X_1}^2,
\ee
by the same argument used to obtain \eqref{en}.
Here, we are using the fact that $\alpha$
satisfies both $|\alpha^{-1}\le C|$ and $|(d/dx)^k\alpha|\le C|\alpha|$
for $k\ge 0$.
Bound \eqref{Nvaren}(i) then follows by \eqref{weightvar} and
the Duhamel formulation
$$
\int \CalN= \int_0^T \big(\int e^{\tilde L(\eps)(T-t)} \partial_x\big)
Q(U,U_x)(t)\, dt,
$$
bounds \eqref{XDuhamel} with $p=\infty$ and $r=2$, and Sobolev embedding
\ba
\|e^{\eta \lx} U(\cdot, t)\|_{L^\infty},
\,
\|e^{\eta \lx} \partial_x U(\cdot, t)\|_{L^\infty}
&\le
\|e^{\eta \lx} U(\cdot, t)\|_{W^{1,\infty}}\\
& \le
C\|e^{\eta \lx} U(\cdot, t)\|_{H^3}.
\ea

Bound \eqref{Nvaren}(ii) follows, similarly, 
by repeating in weighted norm $X_1$
the arguments for \eqref{en} and \eqref{varen} with $U$ replaced by
$W:=U^1-U^2$, where $U^j$ are solutions of \eqref{trunc}
with different initial data $U_0^1$ and $U_0^2$.
%TODO: mention, use Duhamel again?
Bound \eqref{Nvaren}(iii) follows by 
\eqref{Nvaren}(i) together with the equation for $F:=W-e^{\tilde L(\eps)T}W_0$.
%TODO: mention, use Duhamel again, plus lin'd est. t-deriv?
%TODO: a couple more words on Frechet differentiability?
\end{proof}

\br\label{genen}
\textup{
It is readily checked that bounds \eqref{en},\eqref{varen}, \eqref{Nvaren}
%and Proposition \ref{continuity},
hold also under the weakened assumptions described in Remark \ref{gensymm},
in the favorable variables $(U_1, W_2)$,
by rearranging \eqref{sys} as
\ba
\bp U_1\\W_2\ep_t &+ \bp A_{11}(U) & * \\ * & *\ep  \bp U_1\\W_2\ep_x=\\
&\qquad  \qquad
\partial_x \bp 0 & 0 \\ 0 & (\partial W_2/\partial U_2)\beta(U)\ep \bp U_1 \\ W_2\ep_x
+ \bp 0 \\ O(|U_x|^2)\ep
\ea
similarly as in \cite{Z2, Z3, GMWZ} and carrying out energy estimates
as before.
These then imply corresponding estimates in the original 
variables $(U_1, U_2)$, by the observation that,
for $\hat U$ satisfying the linearized equations in $(U_1, U_2)$
coordinates, $(\hat U_1, (\partial W_2/\partial U)(\bar U) \hat U)$
satisfies the linearized equations in $(U_1, W_2)$ coordinates,
and, for $U=\tilde U-\bar U$
satisfying the nonlinear perturbation equation in $(U_1, U_2)$ coordinates,
$(\tilde U_1, W_2(\tilde U))-(\bar U_1, W_2(\bar U))=
( U_1, (\partial W_2/\partial U)(\bar U)  U) + O(|U|^2)$
satisfies the nonlinear perturbation equation in $(U_1, W_2) $ 
coordinates.
}
\er

\section{Bifurcation analysis}\label{bifurcation}

We now carry out the bifurcation analysis following 
the framework of \cite{TZ2}, with a few slight modifications
to simplify the analysis.

\subsection{Construction of the period map}\label{periodmap}

Given a solution $\tilde U$ of \eqref{sys},
define the perturbation variable 
\be\label{pert}
U(x,t,\eps):=\tilde U^\eps (x,t) -\bar U^\eps(x),
\ee
satisfying nonlinear perturbation equations
\be\label{nonlin}
U_t-L_\eps U= Q_\eps (U,U_x)_x,
\qquad
U(x,0,\eps)=U_0(x,\eps).
\ee
By \eqref{trunc}, \eqref{blocktrunc}, 
\be\label{Qbounds}
Q_\eps (U,\eps)=O(|U|^2+|U||U_x|)
%TODO: unneeded, delete:
%\partial_x Q_\eps (U,\eps)=O(|U||U_x|+|U||U_{xx}|+ |U_x|^2)
\ee
so long as $U$ satisfies a uniform $L^\infty$ bound:
in particular, for $\|U\|_{X_1}\le C$.

Decomposing
\be\label{pertdecomp}
U= w_1 \phi_+ + w_2\phi_- + v,
\ee
where $ w_1 \phi_+ + w_2 \phi_- \in \Sigma$,
$v:=\tilde \Pi U\in \tilde \Sigma$,
and coordinatizing as $(w, v)$, 
we obtain after a brief calculation
\ba\label{PHeq2}
\dot w&= \bp \gamma(\eps) & \tau(\eps)\\
-\tau(\eps) & \gamma(\eps)\ep w + N_w(w, v,\eps),\\
\dot v&= \tilde L_\eps v + N_v(w, v,\eps),\\
\ea
where $\tilde L=L\tilde \Pi$ and
\ba\label{Ns}
N_{w,1} \phi_+ + N_{w,2} \phi_- &= \Pi Q_\eps(U,U_x)_x,\\
N_v=\tilde \Pi Q_\eps(U,U_x)_x.
\ea

\begin{lem}\label{Nbds}
$N_w$ is $C^k$, $k\ge 2$, 
from $(w,v)\in \RR^2 \times L^q$
to $\R^2$, for any $0\le q\le \infty$, with
\be\label{quadlin}
|N_v|\le C (|w|^2+\|v\|_{L^q}^2).
\ee
$N_v$ is $C^k$, $k\ge 2$, from 
$(w,v)\in \RR^2 \times X_1$ to $X_2':=
\partial_x \{f:\, \|e^{2\eta \lx} f\|_{H^1}<+\infty\}$ and, 
for $(|w|+\|v\|_{X_1})\le C$, 
$C^1$ from $(w, v)\in \RR^2 \times B_1$
to $\partial_x (L^1)$,
with
\be\label{vquadlin}
\|N_v\|_{X_2'}\le C (|w|^2+\|v\|_{X_1}^2)
\quad \hbox{\rm and }\quad
\|DN_v\|_{\partial_x (L^1)} \le C(|w|+\|v\|_{X_1}). 
\ee
\end{lem}

\begin{proof}
Direct calculation, using \eqref{Qbounds} 
and the $\Pi$-bounds of Lemma \ref{efns}.
\end{proof}

\begin{prop}\label{shortbounds}
For $0\le t\le T$, any fixed $C_1,T>0$, some $C>0$, and $|a|$, $\|b\|_{X_1}$,
$|\eps|$
sufficiently small, system \eqref{PHeq2} with initial data
$(w_0,v_0)=(a, b)$ sufficiently small in $\RR^2\times X_1$ 
 possesses a solution
\be\label{solnop}
(w, v)(a,b,\eps, t)\in \RR^2\times X_1
\ee 
that is $C^{k+1}$ in $t$ and $C^k$ in $(a,b,\eps)$, $k\ge 2$, 
with respect to the weaker norm $B_1$, with
\ba\label{truncshort}
C^{-1}|a|- C\|b\|_{X_1}^2&\le |w(t)|\le C(|a|+ \|b\|_{X_1}^2),\\
\|v(t)\|_{X_1}& \le C(\|b\|_{X_1}+|a|^2),\\
\ea
and
\ba\label{Bshorteq}
|D_{(a,b)}(w, v)(t)|_{\R^2\times B_1
\to B_1}& \le C.
%NOTE: rate determined by linear term $(\|b\|_{X_1}+|a|)$.
\ea

In particular, for $\|b\|_{X_1}\le C_1|a|$, all $0\le t\le T$,
\be\label{truncjust}
\|v(t)\|_{X_1}\le C|w(t)|.
\ee
Likewise, for $\|b\|_{X_1}\le C_1|a|^2$, all $0\le t\le T$,
\be\label{quadbd}
\|v(t)\|_{X_1}\le C|w(t)|^2.
\ee
\end{prop}

\begin{proof}
Existence and uniquess follow by a standard Contraction--mapping
argument, using a priori bounds \eqref{Nvaren}, which
also imply \eqref{truncshort} (by decoupling of linear parts) 
and \eqref{Bshorteq}.
Combining \eqref{truncshort}(i)--(ii),
we obtain evidently \eqref{truncjust} and \eqref{quadbd} 
for $|a|$ sufficiently small.  
\end{proof}

Setting $t=T$ in \eqref{solnop} and applying 
Duhamel's principle/variation of constants,
we may express the period map
\be\label{period}
(a,b,\eps)\to (\hat a, \hat b ):=(w,v)(a,b,\eps, T)
\ee
as a discrete dynamical system 
\ba
\label{PHdyn}
\hat a&=R(\e,T)a+ N_1(a,b,\eps, T),\\
\hat  b&=  S(\eps, T)b+ N_2(a,b,\eps, T)\\
\ea
with $\eps, \, T\in \RR^1$, $a,\, N_1 \in \RR^2$ and $b\in B_1$, $N_2\in B_2$,
where 
\ba\label{RS}
R(\eps, T)&:=e^{\gamma(\eps)T(a,b,\eps)},\\
S(\eps, T)&:=e^{\tilde L_\eps T(a,b,\eps)}\\
\ea
are the linearized solution operators in $w$, $v$
and
\ba\label{PHDuhamel}
N_1(a,b,\eps, T)&:= \int_0^{T} e^{\gamma(\eps)(T-s)}
N_w(w v, \eps)(s)ds,\\
N_2(a,b,\eps, T)&:= \int_0^{T} e^{\tilde L_\eps (T-s)}
N_v(w, v, \eps)(s)ds\\
\ea
the differences between nonlinear and linear solution operators: equivalently,
\ba\label{altN}
N_1(a,b,\eps, T)&= \Pi \CalN(a_1\phi_++a_2\phi_-+b, \eps, T)\\
N_2(a,b,\eps, T)&= (\Id -\tilde \Pi)\CalN(a_1\phi_++a_2\phi_-+b, \eps, T),
\ea
where $\CalN$ is the linearization error $\CalN$ defined in \eqref{N}.

Evidently, periodic solutions of \eqref{PHeq2} 
with period $T$ correspond to fixed points of the 
period map (equilibria of \eqref{PHdyn}) or,
equivalently, zeroes of the displacement map
\ba
\label{displacement}
\Delta_1 (a,b,\eps, T)&:=(R(\e,T)-\Id)a+ N_1(a,b,\eps, T),\\
\Delta_2 (a,b,\eps, T)&:=  (S(\eps, T)-\Id)b+ N_2(a,b,\eps, T).\\
\ea

\subsection{Lyapunov--Schmidt reduction}\label{LS}

We now carry out a nonstandard Lyapunov--Schmidt reduction
following the ``inverse temporal dynamics'' 
framework of \cite{TZ2}, tailored for the situation
that $(S(\eps,T)-\Id)$ is not uniformly invertible, or,
equivalently, $\sigma(\tilde L)$ is not bounded away
from $\{j\pi/T\}$, $j\in \ZZ$.
In the present situation, $\tilde L$ has both an $\ell$-dimensional
kernel (Lemma \ref{ker} below) 
and essential spectra accumulating at $\lambda=0$, and no
other purely imaginary spectra, so that
$(S(\eps,T)-\Id)=(e^{\tilde LT}\tilde \Pi-\Id)$ inherits the same properties;
see \cite{TZ2} for further discusssion.

Our goal, and the central point of the analysis, 
is to solve $\Delta_2(a,b,\eps,T)=0$ for $b$ as
a function of $(a,\eps,T)$, eliminating the transverse variable
and reducing to a standard planar bifurcation problem in the oscillatory
variable $a$.
A ``forward'' temporal dynamics technique would be to rewrite
$\Delta_2=0$ as a fixed point map
\be\label{fwd}
b =  S(\eps, T)b+ N_2(a,b,\eps, T),
\ee
then to substitute for $T$ an arbitrarily large integer multiple $jT$.
In the strictly stable case $\Re \sigma(\tilde L) \le -\eta<0$, 
$|S(\eps, jT)|_{L(X_1, X_1)}<1/2$ for $j$ sufficiently large.
Noting that $N_2$ is quadratic in its dependency, we would have
therefore contractivity of \eqref{fwd} with respect to $b$,
yielding the desired reduction.
However, in the absence of a spectral gap between $\sigma(\tilde L)$
and the imaginary axis, $|S(\eps, jT)|_{L(X_1, X_1)}$ does not
decay, and may be always greater than unity; thus, this naive
approach does not succeed.

The key idea in \cite{TZ2} is to rewrite $\Delta_2=0$ instead 
in ``backward'' form 
\be\label{backwd}
b =  (\Id -S(\eps, T))^{-1} N_2(a,b,\eps, T),
\ee
then show that $(\Id-S)^{-1}$ is well-defined and bounded on $\Range N_2$,
thus obtaining contractivity by quadratic dependence of $N_2$.
Since $(\Id-S)^{-1}N_2$ is formally given by 
$\sum_{j=0}^\infty S^jN_2$
this amounts to establishing convergence:
a stability/cancellation estimate.
Quite similar estimates appear in the nonlinear stability theory,
where the interaction of linearized evolution $S$ and nonlinear
source $N_2$ are likewise crucial for decay.
The formulation \eqref{backwd} can be viewed also as a ``by-hand''
version of the usual proof of the standard Implicit Function Theorem \cite{TZ2}.

\begin{lem}\label{returnsetup}
%TODO: later...
%$R(0,0,\cdot)$ is differentiable for all $\eps$ sufficiently
%small,  
%with $\partial_\eps R(0,0,0)= \partial_\eps \gamma(0) 2\pi /\tau(0)\ne 0$ 
%and $R(0,0,0)=1$; 
Assuming (A1), (A2), (H0)--(H4),  (\D),
$N_1$ is quadratic order and $C^1$ from 
$\RR^2\times L^q\times \RR^2 \to \RR^1$ for any $1\le q\le \infty$,
and $N_2$ is quadratic order from 
$\RR^2\times X_1\times \RR^2 \to X_2$ and $C^1$ from 
$\RR^2\times B_1\times \RR^2 \to B_2$ for $\|b\|_{X_1}$
uniformly bounded, with 
  \begin{equation} \label{ass-N1}
 \begin{aligned}
 |N_1(a,b, \eps, T)|, \quad
 |\partial_{\eps, T}N_1(a,b, \eps, T)|_{L(\R^2, \R^2)}
 & \le C(|a|+\|b\|_{{\cal B}_1})^2, \\
  |\partial_a N_1(a,b, \eps, T)|_{{\cal L}(\R^2, \R^2)} + 
|\partial_b N_1(a,b,\eps, T)|_{{\cal L}(\CalB_1,\RR^2)}
& \le C(|a|+\|b\|_{{\cal B}_1}), \end{aligned}
\end{equation} 
\begin{equation} \label{ass-N2}
 \begin{aligned}\|N_2(a,b, \eps, T)\|_{X_2}, \quad
 |\partial_{\eps, T} N_2(a,b, \eps, T)\|_{L(\R^2, B_2)}
 & \le C(|a|+\|b\|_{X_1})^2,  \\
  |\partial_a N_2(a,b, \eps, T)|_{{\cal L}(\RR^2,\CalB_2)} + |\partial_b 
N_2(a,b, \eps, T)|_{{\cal L}(\CalB_1,\CalB_2)}
& \le C(|a|+\|b\|_{X_1}). \end{aligned} \end{equation} 
\end{lem}

\begin{proof}
Bounds \eqref{ass-N1} follow from representation \eqref{altN},
variational bounds \eqref{Nvaren} of Proposition \ref{Nest}, 
and the $\Pi$-bounds of Lemma \ref{efns}, from which we
likewise obtain
$$
 \begin{aligned}\|N_2(a,b, \eps, T)\|_{X_1} , \quad
 |\partial_{\eps, T} N_2(a,b, \eps, T)\|_{L(\R^2, B_1)}
& \le C(|a|+\|b\|_{X_1})^2,  \\
  |\partial_a N_2(a,b, \eps, T)|_{{\cal L}(\RR^2,\CalB_1)} + |\partial_b 
N_2(a,b, \eps, T)|_{{\cal L}(\CalB_1,\CalB_1)}
& \le C(|a|+\|b\|_{X_1}). \end{aligned} 
$$
The remaining bounds
$$
 \begin{aligned}\|N_2(a,b, \eps, T)\|_{\partial_x e^{-2\eta \lx} H^1 }, \quad
 |\partial_{\eps, T} N_2(a,b, \eps, T)
\|_{L(\R^2, \partial_x e^{-2\eta \lx} H^1) }
& \le C(|a|+\|b\|_{X_1})^2,  \\
  |\partial_a N_2(a,b, \eps, T)|_{{\cal L}(\RR^2,\partial_x L^1)} + |\partial_b 
N_2(a,b, \eps, T)|_{{\cal L}(\CalB_1,\partial_x L^1)}
& \le C(|a|+\|b\|_{X_1}). \end{aligned} 
$$
%TODO: move to later
%From \eqref{Tprop}, we find that $R(0,0,\eps)=e^{2\pi \gamma(\eps)/\tau(\eps)}$
%is $C^k$, $k\ge 2$ for $\eps$ sufficiently small, with $R(0,0,0)=0$
%and $\partial_\eps R(0,0,0)= \partial_\eps \gamma(0)2\pi/\tau(0)\ne 0$
%as claimed.
%NOTE: curiously, NO PROBLEM with this whatsoever...
follow easily from Duhamel representation
\eqref{PHDuhamel}(ii), bounds \eqref{truncshort}(ii), 
\eqref{Bshorteq}, and \eqref{vquadlin},
and bounds \eqref{XDuhamel}, \eqref{varXDuhamel}
on the linearized solution operator $e^{\tilde Lt}\tilde \Pi$.
\end{proof}

 \begin{cor}[\cite{TZ2}] \label{g+} 
Assuming (A1), (A2), (H0)--(H4),  (\D),
 \begin{equation} \label{b} 
\Delta_2(a,b,\eps,T)=0
\qquad (a,b, \eps, T) \in \R^2 \times X_1 \times \RR^2,\end{equation}
 is equivalent to
 \begin{equation} \label{equiv'}
 b = (\Id -S(a,b,\e,T))^{-1}N_2(a,b,\e, T) +\omega
\end{equation} 
for
\begin{equation} \label{omega}
\omega \in 
\kernel(\Id -S(a,b,\e,T)) \cap X_1.
\end{equation} 
\end{cor}
 
 \begin{proof} 
 Applying to the left of \eqref{b} the right inverse $(\Id-S(\e,a,b))^{-1}$ 
given by Proposition \ref{inverse}, we obtain by Lemma \ref{returnsetup}
$$
\begin{aligned}
 \tilde b &:= (\Id -S(\e,a,b))^{-1}(\Id -S(\e,a,b)) b\\
&= (\Id -S(\e,a,b))^{-1}N_2(\e,a,b)\in X_1.
\end{aligned}
$$
Observing that $\tilde b - b$ belongs to $\kernel(\Id -S(\e,a,b)) \cap X_1$
by the right inverse property,
we obtain \eqref{equiv'}. Conversely, \eqref{equiv'} implies 
\eqref{b} by application on the left of $\Id - S(a,b,\e,T).$
\end{proof}

\begin{lem}\label{ker}
The kernel of $(\Id-S(\eps, T))$ is of fixed dimension
$\ell$ as in (H4), is independent of $T$, 
and has a smooth basis $\omega=\bp \omega_1 & \dots & \omega_\ell \ep(\eps)$.
In the Lax case, it is generated entirely by translation invariance.
\end{lem}

\begin{proof}
This follows by the corresponding properties of $\kernel L$ 
assumed in (H4).
\end{proof}

\begin{cor}\label{boundedmap}
%TODO
Assuming (A1), (A2), (H0)--(H4),  (\D), the map
\be\label{CalT}
\CalT(a,b,\eps,T,\alpha)
:=(\Id-S(\eps, T))^{-1}N_2(a,b,\eps,T) + \omega(\eps)\alpha,
\ee
$(\Id-S)^{-1}:X_1 \to X_2$ as defined in Proposition \ref{inverse}
is bounded from $\RR^2\times X_1\times \RR^{2+\ell} \to X_1$ 
and $C^1$ from $\RR^2\times B_1 \times \RR^{2+\ell} \to X_1$ for
$|\alpha|$ bounded and $|a|+\|b\|_{X_1}+|(\eps,T)|$ sufficiently small, 
with
\ba\label{Tbds}
\|\CalT(a,b,\eps,T, \alpha)\|_{X_1}&\le C(|a|+ \|b\|_{X_1}^2),\\
\|\partial_{a,b} \CalT(a,b,\eps,T)\|_{L(B_1, B_1)}&
\le C(|a|+ \|b\|_{X_1}),\\
\|\partial_{T} \CalT(a,b,\eps,T)\|_{L(B_1, B_1)}&
\le C(|a|^2+ \|b\|_{X_1}^2),\\
\|\partial_{\eps} \CalT(a,b,\eps,T)\|_{L(B_1, B_1)}&
\le C(|a|^2+ \|b\|_{X_1}^2 +|\alpha|),\\
\|\partial_{\alpha} \CalT(a,b,\eps,T)\|_{L(B_1, B_1)}&
\le C(|a|^2+ \|b\|_{X_1}^2 +1 ).\\
\ea
\end{cor}

%TODO: old notes, delete later:
%!!!! not $C^1$, with $|a|$ in there... ???
%AHA! Ok, except at $|a|=0$.  But, this will be factored out
%later... giving back $C^1$... or?

\begin{proof}
Immediate, from Proposition \ref{inverse} and Lemmas \ref{returnsetup}
and \ref{ker}.
\end{proof}

\begin{prop}[\cite{TZ2}]\label{LSprop}
Under assumptions (A1)--(A2), (H0)--(H4), (\D) 
(alternatively, the weakened assumptions of Remark \ref{genen}),
there exists a function $\b (a,\eps, T, \alpha)$, bounded
from $\RR^{4+\ell}$ to $X_1$ and $C^1$ from $\RR^{4+\ell}$ to $B_1$,
with
\be\label{Deltasoln}
\Delta_2(a, \b(a,\eps,T,\alpha), \eps, T)\equiv 0,
\ee
\ba\label{betabds}
\|\beta\|_{X_1}, \,
\|\partial_{\eps, T}\beta\|_{L(\R, B_1)}&\le C(|a|^2 + |\alpha|),\\
\|\partial_{a}\beta\|_{L(\R^2, B_1)}&\le C|a|,\\
\|\partial_{\a}\beta\|_{L(\R^2, B_1)}&\le C,\\
\ea
for $|(a,\eps,T,\alpha)|$ sufficiently small.
Moreover, 
for $|(a,\eps,T)|$, $\|b\|_{X_1}$ sufficiently small,
all solutions of \eqref{b} lie on 
the $\ell$-parameter manifold $\{b=\beta(a,\eps,T,\alpha)\}$. 
\end{prop}

\begin {proof}
By Corollary \ref{g+}, \eqref{b} is equivalent to the fixed-point
problem
$$
b= \CalT(a,b,\eps, T,\alpha )
$$
for some $\alpha\in \RR^\ell$.
By \eqref{Tbds}(i)--(ii), 
for $|(a,\eps,T,\alpha)|$, sufficiently small,
$\CalT$ preserves a small ball in $X_1$
on which it is contractive in $b$ with respect to the
weaker norm $\|\cdot\|_{B_1}$.
Observing that closed balls in $X_1$ are closed also in $B_1$,
we may conclude by the contraction-mapping principle the existence
of a unique solution $\beta$, which, moreover, inherits the
regularity of $\CalT$ in its dependence on parameters $(a,\eps,T,\alpha)$.
\end {proof}

\br\label{smoothness}
\textup{
At the expense of further bookkeeping, we may replace \eqref{X1} by
$\|U\|_{X_1}:= \|e^{\eta \lx}U\|_{H^4}$,
$\|\partial_x U\|_{X_2}:=\|\partial_x U\|_{X_1}
+
\|e^{2\eta \lx} U\|_{H^1}$ and carry one further derivative in 
$(a,\eps,T$
throughout the analysis, to obtain $C^2$ dependence of 
$\beta(a,\eps,T,\alpha)$.
Indeed, strengthening (H0) to $k=2r+1$ in (H0),
we may obtain arbitrary smoothness $C^r$ of reduction
function (nullcline) $\beta(\cdot)$.
}
\er

%%%%%%%%%%%%%%%%%%%%%%%%%%%%%%%%%%%%%%%%%%%%%%%%%%%%%%%%%%%%%%%%%%%%%%%%%

%TODO, NOTES:
%NO!  Instead, just do reduction to planar system, for return
%map with $T$ as parameter.  Get condition for reduction map
%$|v|\le C|w|$.  Then quickly sketch planar Poincar\'e,
%just to the point of eliminating $T$. (OR?  What about
%solving for $T$, eliminating $\theta$ equation... better...).

\subsection{Proof of the main theorem}\label{proof}

The bifurcation analysis is straightforward
now that we have reduced to a finite-dimensional problem,
the only tricky point being to deal with the $\ell$-fold 
multiplicity of solutions (parametrized by $\alpha$).
Define to this end
\be\label{tildebeta}
\tilde \beta(a,\eps,T,\hat \alpha):=
\beta(a,\eps,T,|a|\hat \alpha),
\ee
with $\hat \alpha$ restricted to a ball in $\R^\ell$,
noting, by \eqref{betabds}, that
\ba\label{tildebetabds}
\|\tilde\beta\|_{X_1}, \,
\|\partial_{a,\eps, T, \ha}\tilde\beta\|_{L(\R, B_1)}&\le C|a|,
\ea
with $\tilde \beta$ Lipshitz in $(a,\eps,T,\ha)$ and $C^1$ away from $a=0$. 
Solutions $(w,v)$ of \eqref{PHeq2} originating at $(a,b)=(a,\tb)$,
by \eqref{truncjust}, remain for $0\le t\le T$ in a cone
\be\label{Cdef}
\CalC:=\{(w,v):\, |v|\le C_1 |w|\},
\ee
$C_1>0$.
% fixed.
Likewise, any periodic solution of \eqref{PHeq2} originating
in $\CalC$, since it necessarily satisfyies $\Delta_2=0$, 
must originate from data $(a,b)$ of this form.

\begin{proof}[Proof of Theorem \ref{main}]
Defining $b\equiv \tilde \beta(a,\eps,T,\ha)$, and recalling
invariance of $\CalC$ under flow \eqref{PHeq2}, we may view
$v(t)$ as a multiple
\be\label{multiple}
v(t)=\gamma(a,\eps,T,\ha, t) w(t)
\ee
of $w(t)$, where $\gamma$ is bounded, Lipshitz in all arguments,
and $C^1$ away from $a=0$.
Substituting into \eqref{PHeq2}(i), we obtain a planar ODE 
\be\label{nonautonomous}
\dot w = \bp \gamma(\eps) & \tau(\eps)\\
-\tau(\eps) & \gamma(\eps)\ep w + \tilde N(w,\eps, T, t, \ha, a)
\ee
in approximate Hopf normal form, with nonlinearity 
$\tilde N:=N_w(w, v, \eps )$ 
now nonautonomous and depending on the additional parameters $(T,\ha,a)$, but 
still satisfying the key bounds
\ba\label{tildeNbds}
|\tilde N|, \,
|\partial_{\eps, T, \ha, a} \tilde N|&\le C|w|^2;
\quad
|\partial_{w} \tilde N|&\le C|w|
\ea
along with planar bifurcation criterion (\D)(ii).
From \eqref{tildeNbds}, we find that $\tilde N$ is $C^1$ in all
arguments, also at $a=0$.
%NOTE: HERE, we encounter and repair Lip. trouble for first time...
%We'll need to do it again for $r-\theta$ coord's...
%OR?  Could we do return map a different way maybe??? (in std. coord's,
%thus avoiding introduction of singularities???)
%NOTE, also: this inherently limits us to $C^1$ in doing whole family.
%BUT: not for single soln. $\ha\equiv0$...
By standard arguments (see, e.g., \cite{HK, TZ1}), 
we thus obtain a classical Hopf bifurcation in the variable $w$
with regularity $C^1$, yielding existence and uniqueness up to
time-translates of an $\ell$-parameter family of solutions 
originating in $\CalC$, indexed by $r$ and $\delta$
with $r:=a_1$ and (without loss of generality) $a_2 \equiv 0$.
%except for uniqueness up to spatial translates
%in the case of a Lax-type background shock,
%with $r:=a_1$ and (without loss of generality) $a_2 \equiv 0$.
%In particular, this completes the proof in the overcompressive case.
It remains only to establish uniqueness up to spatial translates.

In the Lax case, we observe, first, that, by dimensional considerations,
the one-parameter family constructed must agree with the one-parameter
family of spatial translates.  Second, we argue as in \cite{TZ2}
%TODO: precise ref. here
that any periodic solution has a spatial translate originating
in $\CalC$, yielding uniqueness up to translation among all solutions
and not only those originating in $\CalC$; 
%see Section 2.21, \cite{TZ2} for further details.
see Proposition 2.17 and Corollary 2.18 \cite{TZ2} for further details.
%TODO: recheck ref.

In the overcompressive case, we observe, likewise, that the $\ell$-parameter
family constructed consists of spatial translates of a smooth
$(\ell-1)$-parameter family of distinct orbits of traveling-wave
ODE \eqref{ode}, then argue as in \cite{TZ2}
%TODO: precise ref. here
that any periodic solution has a spatial translate lying within
a corresponding cone $\CalC'$ about some member of this 
$(\ell-1)$-parameter family.
Constructing solutions about each such member by the same technique, 
we thus obtain from their union an $(\ell-1)$ parameter family containing
translates of all periodic solutions on an 
entire neighborhood of $\bar U$, and not only in cone $\CalC$.
%See \cite{TZ2} for details.
%TODO: precise ref. here
See Proposition 2.22 and Corollary 2.23 \cite{TZ2} for further details.
\end{proof}

\begin{rems}\label{regrems}
\textup{
The apparent restriction to $C^1$ regularity caused by factor
$|a|$ in the definition of $\tilde \b$ is illusory, since
we restrict eventually to the ray $a_2\equiv 0$ and $a_1\ge 0$,
on which $\tb$ is as smooth as $\b$. 
By Remark \ref{smoothness}, the latter may be made as smooth
as desired by assuming sufficient regularity in (H0), and so
we can carry out a bifurcation analysis to arbitrary desired regularity.
}

\textup{
2. 
%By a standard cancellation argument \cite{HK}, 
%it may be shown that
%$\|N_1\|_{X_1}=O(|a|^3+\|b\|_{X_1}^2)$.
Restricting to the ``central'' solutions $\ha\equiv 0$
of the constructed 
%$\ell$-dimensional manifold 
cone of solutions
yields $\tb=O(|a|^2)$, and thus the exact Hopf normal form
\be\label{refnormal}
\dot w = \bp \gamma(\eps) & \tau(\eps)\\
-\tau(\eps) & \gamma(\eps)\ep w + N_w(w,0,\eps)+ \tilde M(w,\eps, T, t, \ha, a),
\ee
%in place of of \eqref{nonautonomous},
where $\tilde M=O(|w|^3)$ and $\partial_{w,a}M=O(|w|^2)$.
For these central solutions,
we thus recover all of the standard description of Hopf bifurcation:
in particular, that 
$\eps(r)$ is a pitchfork bifurcation with $(d\eps/dr)(0)=0$
(hence $\sgn (d^2\eps/dr^2)(0)$ is expected to determine stability;
see the final paragraph of Section \ref{discuss}).
%TODO: how though to deal with weirdness of polar coordinates, even
%with standard, planar Hopf???  Not so clear... (though now I do at
%least see $C^2$...).
%That is, for the ``central'' solution $\ha\equiv 0$,
%of the constructed $\ell$-dimensional manifold of solutions, 
%we recover all of the standard description of Hopf bifurcation.
}
\end{rems}

\section{Physical applications}\label{applications}

We conclude by describing applications to various 
physical systems.
To generate example systems and waves, it is easier to look
for traveling-wave solutions $U(x,t)=\bar U(x-st)$
with possibly nonzero speed, to be transformed to steady solutions
by the change of coordinates $x\to x-st$ as described in the introduction.
Examples satisfying (A1)--(A2), (H0)--(H4) are: 
\medskip

1. The general Navier--Stokes equations of compressible gas dynamics, written 
in Lagrangian coordinates, appear as
\begin{equation}
\left\{\begin{array}{l}
 v_t -u_x = 0,\\
 u_t + p_x =((\nu/v) u_{x})_x,\\
 (e+u^2/2)_t + (pu)_x = ((\kappa/v) T_x + (\mu/v) uu_x)_x,
\end{array}\right.
\label{NS}
\end{equation}
where $v>0$ denotes specific volume, $u$ velocity, $e>0$ internal energy, 
$T=T(v,e)>0$ temperature, 
%$T=c^{-1}e>0$ temperature, 
$p=p(v,e)$ pressure, and $\mu>0$ and $\kappa>0$
are coefficients of viscosity and heat conduction, respectively.
For simplicity, assume an ideal temperature dependence
\be\label{ideal}
T=T(e),
\ee
independent of $v$.  This assumption
can be removed with further effort; see Remark \ref{superfluous}.
%  $de=T ds -pdv$
Defining $U_1:=(v)$, $U_2:=(u,e+|u|^2/2)$, we find that
the single hyperbolic mode is $v$, and the associated equation is
linear by inspection;  likewise, the viscosity matrix 
$b=\bp  \nu/v & 0\\
* & \kappa/\nu c\ep$
is lower triangular with positive real diagonal entries, 
hence its spectrum is positive real and (by Lyapunov's Lemma)
there exists
$A^0_{22}$ symmetric positive definite such that $A^0_{22}b>0$.
Conditions (A1)--(A2) and (H1)--(H4)
%TODO: check H numbers
are thus satisfied under the 
mild assumptions of monotone temperature-dependence,
\be\label{Tstab}
T_e>0,
\ee
and thermodynamic stability of the endstates, 
\be\label{thermstab}
(p_v)_\pm <0,  (T_e)_\pm >0;
\ee
see, e.g., \cite{MaZ4, Z2, Z3} for further discussion.
Notably, this allows the interesting case of a van der Waals-type 
equation of state, with $p_v>0$ for some values of $v$ along
the connecting profile.

2. Next, consider the equations of MHD:
\begin{equation}
\left\{\begin{array}{l}
v_t -u_{1x} = 0,\\
 u_{1t} + (p+ (1/2\mu_0)(B_2^2+B_3^2))_x =((\nu/v) u_{1x})_x,\\
 u_{2t}  - ((1/\mu_0)B_1^*B_2)_x =((\nu/v) u_{2x})_x,\\
 u_{3t}  - ((1/\mu_0)B_1^*B_3)_x =((\nu/v) u_{3x})_x,\\
 (vB_2)_{t}  - (B_1^*u_2)_x =((1/\sigma\mu_0 v) B_{2x})_x,\\
 (vB_3)_{t}  - (B_1^*u_3)_x =((1/\sigma\mu_0 v) B_{3x})_x,\\
(e+(1/2)(u_1^2+ u_2^2+u_3^2) + (1/2\mu_0)v(B_2^2+B_3^2))_t \\
 \qquad\qquad + [(p+ (1/2\mu_0)(B_2^2+B_3^2))u_1 - (1/\mu_0)B_1^*
  (B_2u_2+B_3u_3)]_x \\
 \qquad\quad
  = [(\nu/v)u_1u_{1x} + (\mu/v) (u_2u_{2x} + u_3 u_{3x})\\
 \qquad\qquad +(\kappa/v) T_x + 
  (1/\sigma \mu_0^2 v)(B_2B_{2x} + B_3 B_{3x})]_x,
\end{array}\right.
\label{MHD}
\end{equation}
where $v$ denotes specific volume, $u=(u_1,u_2,u_3)$ velocity,
$p=P(v,e)$ pressure, $B=(B_1^*,B_2,B_3)$ magnetic induction,
$B_1^*$ constant, $e$ internal energy, $T=T(v,e)>0$ temperature, and
$\mu>0$ and $\nu>0$ the two coefficients of viscosity, $\kappa>0$
the coefficient of heat conduction, $\mu_0>0$ the magnetic
permeability, and $\sigma>0$ the electrical resistivity.
Under assumptions \eqref{ideal}, \eqref{Tstab}, and \eqref{thermstab},
conditions (A1)--(A2) are again satisfied, and
conditions (H1)--(H4) are satisfied (see \cite{MaZ4}) 
under the generically satisfied
assumptions that the shock be of Lax or overcompressive type, 
the endstates $U^\eps_\pm$
be strictly hyperbolic, and the speed $s$ be nonzero, 
i.e., the shock move with nonzero speed relative to the background 
fluid velocity, with $U_1:=(v)$, 
$U_2:=
(u, B, e+|u|^2/2 + v|B|^2/2\mu_0)$.
(For gas dynamics, only Lax-type shocks and nonzero speeds can occur,
and all points $U$ are strictly hyperbolic.)

3. (MHD with infinite resistivity/permeability)
An interesting variation of (\ref{MHD}) that is of interest in certain 
astrophysical parameter regimes is the limit in which either electrical 
resistivity $\sigma$, magnetic permeability $\mu_0$, or both, go to infinity,
in which case the righthand sides of the fifth and sixth equations
of (\ref{MHD}) go to zero and there is a three-dimensional
set of hyperbolic modes $(v,vB_2,vB_3)$ instead of the usual one.
By inspection, the associated equations are still linear in the
conservative variables.
Likewise, (A1)--(A2), (H1)--(H4) hold 
under assumptions \eqref{ideal}, \eqref{Tstab}, and \eqref{thermstab}
for nonzero speed Lax- or overcompressive-type shocks 
with strictly hyperbolic endstates. 

4. (multi-species gas dynamics or MHD)
Another simple example for which the hyperbolic modes are vectorial
is the case of miscible, multi-species flow, neglecting species diffusion,
in either gas- or magnetohydrodynamics.
In this case, the hyperbolic modes consist of $k$ copies of the hyperbolic
modes for a single species, where $k$ is the number of total species.
Again, the associated equations are linear, and 
(A1)--(A2), (H1)--(H4) hold for nonzero speed Lax- or overcompressive-type
shocks with strictly hyperbolic endstates, 
under assumptions \eqref{ideal}, \eqref{Tstab}, and \eqref{thermstab}.

\br\label{superfluous}
\textup{
Assumption \eqref{ideal} may be removed in all cases by
working in the more general framework described in Remarks 
\ref{gensymm} and \ref{genen},  
with variables $U_1=(v)$, $W_2=(u, T)$ (resp. $W_2=(u,B,T)$).
%See \cite{Z2, Z3} for similar analysis in Eulerian coordinates.
}
\er

%%%%%%%%%%%%%%%%%%%

%APPENDIX

%%%%%%%%%%%%%%
%
\appendix

\section{Lagrangian vs. Eulerian formulation}\label{eulerian}
In this appendix, we discuss a bit further the 
role of structural assumptions (A1)--(A2) in the
energy estimates of Section \ref{energy}.
Estimate \eqref{en} holds by a similar, quasilinear version of the
argument of Proposition \ref{energyest}, 
for the general class of hyperbolic--parabolic systems
considered in \cite{Z3}.
However, \eqref{varen} is much more delicate,
as are the related estimates
\eqref{en},\eqref{varen}, and \eqref{Nvaren}.

To see why, note that \eqref{varen} is essentially a variational bound,
measuring the difference between two solutions, and variational
bounds typically cost an additional derivative.
For example, in the quasilinear symmetric hyperbolic case,
we have $\|U(t)\|_{H^s}\le C\|U_0\|_{H^s}$ by standard energy
estimates, whereas the same energy estimates applied to 
the difference between two solutions $U^1$ and $U^2$
yield, rather,
$$
\|U^1-U^2\|_{H^s}(t)\le C\|U^1_0-U^2_0\|_{H^s}
\|U^1_0\|_{H^{s+1}}.
$$
In the linear case, 
$\|U^1-U^2\|_{H^s}(t)\le C\|U^1_0-U^2_0\|_{H^s}^2$ follows by
superposition, and so this issue does not arise; this motivates
the assumption of linearity of the hyperbolic part of \eqref{sys}
made in (A1).

Indeed, when (A1) fails, so, typically, does \eqref{varen}.
For example, consider the Eulerian version \eqref{Eulereg}
of the isentropic compressible Navier--Stokes equations,
rewritten in quasilinear form
\begin{equation}
\left\{\begin{array}{l}
 \rho_t + u\rho_x + \rho u_x=0,\\ 
u_t + u u_x + \rho^{-1}p'(\rho)\rho_x =\nu\rho^{-1} u_{xx}.
\end{array}\right.
\label{symmEulereg}
\end{equation}
Perturbing about the constant solution $(\rho, u)\equiv (1,0)$,
and assuming without loss of generality that $p'(1)=1$,
we obtain perturbation equations
\begin{equation}
\left\{\begin{array}{l}
 \rho_t + u\rho_x + u_x= -\rho u_x,\\
 u_t + \rho_x - \nu u_{xx}=O(|\rho|^2 + |u|^2+ |\rho| |u_x|)_x
\end{array}\right.
\label{pertEulereg}
\end{equation}
so long as $\|\rho\|_{L^\infty}\le C\|\rho\|_{H^1}$
remains sufficiently small.
Applying the same energy estimates as in the proof of Proposition
\ref{energyest}, we easily obtain 
\be\label{Euleren}
\|(\rho, u)(\cdot, T)\|_{H^s}^2+\int_0^T \|u(\cdot, t)\|_{H^{s+1}}^2\, dt
\le C\|(\rho_0, u_0)\|_{H^s}^2\\
\ee
for $s\ge 1$, similarly as in \eqref{en}.

Now, consider the solution $\hat \rho$ of
$$
\hat \rho_t + u\hat \rho_x +u_x =0,
$$
with $\rho$, $u$ as determined by \eqref{pertEulereg}.
Defining $\hat E:=\hat \rho -\rho$, we have
$$
\hat E_t + u\hat E_x= \rho u_x.
$$
Applying the same energy estimates carried out on the variation $E$
in the proof of Proposition \ref{energyest}, we thus obtain
for $s\ge 1$ that
\ba
\|\hat \rho - \rho\|_{H^s}(t)&\le 
C\|(\rho_0, u_0)\|_{H^s}^2. 
%+\int_0^t \|u\|_{H^{s+1}}(r)dr\\
%&\le
%C\|(\rho_0, u_0)\|_{H^s}^2 + t\int_0^t \|u\|_{H^{s+1}}^2(r)dr\\
\ea

On the other hand, defining $\tilde \rho$ to be the solution
of the linearized equation
$$
\tilde \rho_t + u_x=0,
$$
%and $\tilde E:=\tilde \rho - \hat \rho$, 
we find
\ba\label{one}
\|\tilde \rho(t) - \rho_0\|_{H^s} \le \int_0^t \|u\|_{H^{s+1}}(r)dr
&\le t^{1/2} \Big(\int_0 \|u\|_{H^{s+1}}^2 (r)dr\Big)^{1/2}\\
&\le Ct^{1/2} \|(\rho_0, u_0)\|_{H^s}.
\ea
Likewise, defining $\bar \rho$ to be the solution of
$$
\bar \rho_t + u\bar \rho_x=0,
$$
and $\bar e:=\bar \rho -\hat \rho$
we find from similar energy estimates, together with
$$
\partial_t \|\partial_x^\ell \bar e\|_{L^2}=
\partial_t \|\partial_x^\ell \bar e\|_{L^2}^2/
2\|\partial_x^\ell \bar e\|_{L^2},
$$
that
\be\label{two}
\|\bar \rho -\hat \rho\|_{H^s}(t)\le \int_0^t \|u\|_{H^{s+1}}(r)dr
\le Ct^{1/2}\|(\rho_0, u_0)\|_{H^s}.
\ee

Combining \eqref{Euleren}, \eqref{one}, and \eqref{two}, we
obtain
$$
\|\rho-\tilde \rho\|_{H^s}(t)=
\|\bar \rho-\rho_0\|_{H^s}(t) + O(t^{1/2})\|(\rho_0,u_0)\|_{H^s},
$$
which, for $t$ sufficiently small, implies that
$\|\rho-\tilde \rho\|_{H^s}(t)\le C \|(\rho_0,u_0)\|_{H^s}^2$ only if
$\|\bar \rho-\rho_0\|_{H^s}(t)\le C \|(\rho_0,u_0)\|_{H^s}^2$.
But, $\bar \rho$, as the solution to a simple transport equation, 
is exactly a displacement 
$$
\bar \rho(x,t)= \rho_0(X(x,t))
$$
of $\rho_0$,
where $X$ satisfies characteristic equation $\partial_t X= u(X)$,
and so in the absence of derivative bounds on $\partial_x^s \rho_0$, 
$\|\bar \rho(t) -\rho_0\|_{H^s}$ 
is in general no smaller than $2\|\rho_0\|_{H^s}$.
Thus, $\|\rho-\tilde \rho\|_{H^s}(t)\le C \|(\rho_0,u_0)\|_{H^s}^2$ 
in general does not hold in the Eulerian formulation.

%NOTE: don't need nonzero speed, or any of finer details.
%This estimate is by comparison quite crude!  The business
%with transport, however, is interesting and subtle.

\medskip
The resolution of this apparent paradox is that the Lagrangian
formulation incorporates nonlinear transport effects into the
equations via a change of independent variables, so that
$\|\cdot\|_{H^s}$ compares different solutions along their
respective particle-paths rather than at a fixed location,
thus eliminating the principal contribution computed above.

\br\label{barrier}
\textup{
Some such superlinear bound as
\eqref{quaderror} appears necessary for bifurcation analysis
based on temporal dynamics,
at least as usually performed based on autonomous linearized equations.
This is an obstruction to the application of 
invariant manifold/bifurcation techniques to hyperbolic systems 
such as the incompressible Euler equations \cite{Li}.
%TODO: BUT! check new Schlag paper.  Apparently they DO carry
%out such analysis... ??? interesting... add ref. if appropriate.
Likewise, for hyperbolic--parabolic systems, 
it limits us at present to the Lagrangian
formulation and a single spatial dimension.
For related discussion, see \cite{Ho1, Ho2, HoZ, D}.
}
\er

\end{document}